\documentclass{amsart}

\usepackage{amssymb}

\usepackage{mathtools}
\mathtoolsset{showonlyrefs}

\numberwithin{equation}{section}

\usepackage{enumerate}
\usepackage{enumitem}
\usepackage{refcount}

\newenvironment{eqenumerate}
{\begin{enumerate}[ref=\thesection.\theenumi]
		
		\setcounter{enumi}{\value{equation}}}
	{\setcounter{equation}{\value{enumi}}
\end{enumerate}}

\setlist[itemize]{align=parleft,left=1em..2em}

\usepackage{aligned-overset}
\usepackage{mleftright}

\theoremstyle{plain}

\newtheorem{theorem}{Theorem}[section]
\newtheorem{proposition}[theorem]{Proposition}
\newtheorem{quotetheorem}{Theorem}[section]

\newtheorem{lemma}[theorem]{Lemma}
\newtheorem{corollary}[theorem]{Corollary}

\theoremstyle{definition}

\newtheorem{definition}[quotetheorem]{Definition}
\newtheorem{problem}{Problem}

\theoremstyle{remark}

\newtheorem{remark}[theorem]{Remark}

\newenvironment{proofof}[1]
{\begin{proof}[Proof of #1]}
	{\end{proof}}

\newcommand{\set}[2]{\{ \, #1 : #2 \, \} } 

\newcommand{\C}{\mathbb{C}} 
\newcommand{\R}{\mathbb{R}} 
\newcommand{\Z}{\mathbb{Z}} 
\newcommand{\N}{\mathbb{N}} 

\newcommand{\SO}{\text{SO}} 

\newcommand{\hilbert}{L^2(\R^d)}

\NewDocumentCommand\sobolev{ m }{ 
	\ifx#10%
	H^{s_0 / 2}(\R^d)%
	\else%
	{\ifx#11%
		H^{s/2}(\R^d)%
		\else%
		{\ifx#12%
			H^{s}(\R^d)%
			\else%
			H^{#1}(\R^d)
			\fi}
		\fi}
	\fi}

\NewDocumentCommand\Linfty{ D(){\R^d} }{ L^\infty(#1) }

\NewDocumentCommand\fracLaplacian{ O{1} o{s} }{
	\ifx#11%
	( - \Laplacian + 1)^{s/4}%
	\else
	{
		\ifx#12%
		( - \Laplacian + 1)^{s/2}%
		\else%
		{\ifx#1-%
			( - \Laplacian + 1)^{-s/4}%
			\else
			{
				\ifx#10%
				( - \Laplacian + 1)^{s_0/4}%
				\else
				( - \Laplacian + 1)^{#1}%
				\fi}
			\fi}
		\fi}
	\fi}

\DeclareMathOperator*{\essinf}{ess\,inf} 
\DeclareMathOperator{\sgn}{sgn} 
\DeclareMathOperator{\supp}{supp} 

\NewDocumentCommand\Fourier{ o }{ 
	{
		\IfNoValueF{#1}{\mathcal{F}{#1}}
		\IfNoValueT{#1}{\mathcal{F}}
	}
}

\newcommand{\abs}[1]{\lvert#1\rvert} 
\newcommand{\norm}[1]{\lVert#1\rVert} 
\newcommand{\Norm}[1]{\mleft \lVert#1 \, \mright \rVert} 
\newcommand{\innerproduct}[1]{\langle #1 \rangle}

\newcommand{\indicator}[1]{\textbf{1}_{#1}}

\newcommand{\measure}[1]{ \lvert#1\rvert }

\renewcommand{\Re}{\operatorname{Re}} 

\NewDocumentCommand\adjoint{m}{ #1^* }

\newcommand{\Laplacian}{ { \Delta } }

\newcommand{\dampco}{b}
\newcommand{\dampop}{B}

\newcommand{\levelofdampco}{\beta}
\NewDocumentCommand\setS{ O{\dampco} O{\levelofdampco} }{ S(#1, #2) }
\newcommand{\dampconorm}{\norm{\dampco}_{\Linfty}}
\newcommand{\dampopnorm}{\norm{\dampop}_{\hilbert \to \hilbert}}

\newcommand{\levelofm}{\varepsilon}
\NewDocumentCommand\setsigma{ O{m} O{\levelofm} }{ \Sigma(#1, #2) }
\NewDocumentCommand\setSigma{ O{\lambda} O{s} O{\levelofm} }{ \setsigma[#1, #2][#3] }

\NewDocumentCommand\gen{ O{\dampop} O{s} }{{A}(#2, #1)}

\NewDocumentCommand\semigroup{ O{\gen} }{ e^{ t #1} }

\NewDocumentCommand\anihiconst{ O{S} O{\setSigma} }{ \textbf{\upshape C}_{\text{\upshape Anni.} }(#1, #2) }

\NewDocumentCommand\diagonizer{ m O{\lambda} O{s} }{ P_{#1}(#2, #3) }
\NewDocumentCommand\orderq{ O{s} O{s_0} O{p} }{q(#1, #2, #3)}
\newcommand{\modulation}[1]{ M_{#1} }

\usepackage{hyperref}
\hypersetup{
	setpagesize=false,
	bookmarksnumbered=true,
	bookmarksopen=true,
	colorlinks=true,
	linkcolor=blue,
	citecolor=blue,
}

\usepackage[abbrev, alphabetic, y2k, msc-links]{amsrefs}

\BibSpec{arXiv}{%
	+{}{\PrintAuthors}{author}
	+{,}{ \textit}{title}
	+{,}{ arXiv: }{eprint}
	+{}{ \parenthesize}{date}
	+{.}{}{transition}
}

\begin{document}
	
	\title[Uncertainty principle and damped Klein--Gordon equation]{The uncertainty principle and energy decay estimates of the fractional Klein--Gordon equation with space-dependent damping}
	
	\author{Soichiro Suzuki}
	\address{Graduate School of Mathematics, Nagoya University, Furocho, Chikusaku, Nagoya, Aichi, 464-8602, Japan}
	\email{m18020a@math.nagoya-u.ac.jp}
	\thanks{This work was supported by Japan Society for the Promotion of Science (JSPS) KAKENHI Grant Number 20J21771.}
	\subjclass[2020]{42B37}
	
	\begin{abstract}
		We consider the $s$-fractional Klein--Gordon equation with space-dependent damping on $\mathbb{R}^d$. 
		Recent studies reveal that the so-called geometric control conditions (GCC) are closely related to semigroup estimates of the equation. 
		Particularly, in the case $d = 1$, a necessary and sufficient condition for the exponential stability in terms of GCC is known for any $s > 0$.
		On the other hand, in the case $d \geq 2$ and $s \geq 2$, Green--Jaye--Mitkovski (2022) proved that an `$1$-GCC' is sufficient for the exponential stability, 
		but also conjectured that it is not necessary if $s$ is sufficiently large.
		In this paper, we prove the equivalence between the exponential stability and a kind of the uncertainty principle in Fourier analysis.
		As a consequence of the equivalence, we show that the $1$-GCC is not necessary for the exponential stability in the case $s \geq 4$. 
		Furthermore, we also establish an extrapolation result with respect to $s$.
		In particular, we can obtain the polynomial stability for the non-fractional case $s = 2$ from the exponential stability for some $s > 2$. 
	\end{abstract}

	\maketitle
	
	\section{Introduction}
	We consider the following damped fractional Klein--Gordon equation introduced by \cite{MR4063985}:
	\begin{equation} \label{eq:DFKGeq}
		\begin{cases}
			u_{tt} (x, t) + \dampco(x) u_t(x, t) + \fracLaplacian[2] u(x, t) = 0 , & (x, t) \in \R^d \times \R_{> 0} , \\
			u(x, 0) = u_0(x) , & x \in \R^d , \\
			u_t(x, 0) = v_0(x) , & x \in \R^d , 
		\end{cases}
	\end{equation}
	where $s > 0$, $0 \leq \dampco \in \Linfty$ and $(u_0, v_0) \in \sobolev{1} \times \hilbert$.
	The energy of the solution of \eqref{eq:DFKGeq} at time $t \geq 0$ is defined by 
	\begin{align}
		E(t) 
		&\coloneqq \norm{ ( u(\,\cdot\,, t) , u_t(\,\cdot\,, t) ) }_{\sobolev{1} \times \hilbert}^2 \\
		&= \int_{x \in \R^d} ( \abs{ ( \fracLaplacian[1] u(x, t)  }^2 + \abs{ u_t(x, t)}^2 ) \, dx .
	\end{align}
	It is easy to check that
	\begin{equation}
		\frac{d}{dt} E(t) = -2 \int_{x \in \R^d} \dampco(x) \abs{u_t(x, t)}^2 \, dx \leq 0
	\end{equation}
	holds (at least formally) and thus the energy $E(t)$ is non-increasing.
	In particular, if $\dampco(x) \equiv 0$, i.e.\ there is no damping, then the energy is conserved.
	Also, if $\dampco(x) \equiv \dampco_0 > 0$, then we can show that
	\begin{equation} \label{eq:energy estimate for positive constant damping}
		E(t) \lesssim (1 + t)^{2} e^{- \dampco_0 t} E(0)
	\end{equation}
	holds. 
	In fact, the solution of \eqref{eq:DFKGeq} can be written explicitly via the Fourier transform $\Fourier$:
	\begin{align} \label{eq:explicit solution of DFKG}
		\Fourier[u](\xi, t) 
		&= e^{- \dampco_0 t / 2} \mleft( \dampco_0 \cosh{( \varphi(\xi) t )} + \frac{t}{2} \frac{\sinh{( \varphi(\xi) t )} }{ \varphi(\xi) t } \mright) \Fourier[u_0](\xi) \\
		&\quad + t e^{- \dampco_0 t / 2} \frac{ \sinh{( \varphi(\xi) t )} }{ \varphi(\xi) t } \Fourier[v_0](\xi) ,
	\end{align}
	where 
	\begin{equation}
		\varphi(\xi) \coloneqq \frac{ \sqrt{ \dampco_0^2 - 4 ( \abs{\xi}^2 + 1 )^{s/2} } }{2} ,
	\end{equation}
	and \eqref{eq:explicit solution of DFKG} implies \eqref{eq:energy estimate for positive constant damping}.
	Obviously this method does not work in the non-constant case, thus we need to find another approach.
	Let define $\gen[\dampco] \colon \sobolev{2} \times \sobolev{1} \to \sobolev{1} \times \hilbert$ by 
	\begin{equation}
		\gen[\dampco] \coloneqq \begin{pmatrix}
			0 & 1 \\
			- \fracLaplacian[2] & - \dampco
		\end{pmatrix} .
	\end{equation}
	Then it generates a contraction semigroup $\{ \semigroup[\gen[\dampco]] \}_{t \geq 0}$ on $\sobolev{1} \times \hilbert$ and the solution and energy of \eqref{eq:DFKGeq} are given by
	\begin{equation}
		(u, u_t)
		= \semigroup[\gen[\dampco]] (u_0, v_0) , \quad 
		E(t) = \norm{ \semigroup[\gen[\dampco]] (u_0, v_0) }_{\sobolev{1} \times \hilbert}^2 .
	\end{equation}
	Therefore, we will investigate the operator norm of $\{ \semigroup[\gen[\dampco]] \}_{t \geq 0}$. 
	The following equivalences between semigroup estimates and resolvent estimates are immediate from semigroup theory.
	\begin{quotetheorem}[\cite{MR834231}*{Theorem 3}, \cite{MR4215997}*{Theorem 1.3}] \label{theorem:exp stable = uniform resolvent of A}
		Let $s > 0$ and $0 \leq \dampco \in \Linfty$. 
		Then the following are equivalent:
		\begin{eqenumerate}
			\item \label{item:exp stable}
			$\{ \semigroup[\gen[\dampco]] \}_{t \geq 0}$ is exponentially stable, that is, there exist $M,  \omega > 0$ such that
			\begin{equation}
				\norm{ \semigroup[\gen[\dampco]] }_{ \sobolev{1} \times \hilbert \to \sobolev{1} \times \hilbert } 
				\leq M e^{ - \omega t} 
			\end{equation}
			holds for any $t \geq 0$.
			\item \label{item:uniform resolvent of A}
			There exists $C > 0$ such that
			\begin{equation} 
				\norm{F}_{\sobolev{1} \times \hilbert} 
				\leq C \norm{ ( \gen[\dampco] - i \lambda ) F }_{\sobolev{1} \times \hilbert}
			\end{equation}
			holds for any $\lambda \in \R$ and $F \in \sobolev{2} \times \sobolev{1}$. 
		\end{eqenumerate}
		Furthermore, we have the following:
		\begin{itemize}
			\item 
			If \eqref{item:exp stable} holds with $(M, \omega) = (M_0, \omega_0)$, 
			then \eqref{item:uniform resolvent of A} holds with $C = M_0 \omega_0^{-1}$.
			\item 
			If \eqref{item:uniform resolvent of A} holds with $C = C_0$, 
			then \eqref{item:exp stable} holds with $(M, \omega) = ( e^{\pi / 2} , C_0^{-1} )$.
		\end{itemize}
	\end{quotetheorem}
	\begin{quotetheorem}[\cite{MR2606945}*{Theorem 2.4}] \label{theorem:poly stable = poly resolvent of A}
		Let $s, p > 0$ and $0 \leq \dampco \in \Linfty$. 
		Then the following are equivalent:
		\begin{eqenumerate}
			\item \label{item:poly stable}
			$\{ \semigroup[\gen[\dampco]] \}_{t \geq 0}$ is $1 / p$-polynomially stable, that is, there exists $M > 0$ such that
			\begin{equation}
				\norm{ \semigroup[\gen[\dampco]] \gen[\dampco]^{-1} }_{ \sobolev{1} \times \hilbert \to \sobolev{1} \times \hilbert } 
				\leq M ( 1 + t )^{- 1 / p} 
			\end{equation}
			holds for any $t \geq 0$.
			\item \label{item:poly resolvent of A}
			There exists $C > 0$ such that
			\begin{equation} 
				\norm{F}_{\sobolev{1} \times \hilbert} 
				\leq C (1 + \abs{\lambda})^{p} \norm{ ( \gen[\dampco] - i \lambda ) F }_{\sobolev{1} \times \hilbert}
			\end{equation}
			holds for any $\lambda \in \R$ and $F \in \sobolev{2} \times \sobolev{1}$. 
		\end{eqenumerate}
	\end{quotetheorem}
	\begin{quotetheorem}[\cite{MR1618254}*{TH\'{E}OR\`{E}ME 3}, \cite{MR2460938}*{Proposition 1.3}] \label{theorem:log stable = exp resolvent of A}
		Let $s, p > 0$ and $0 \leq \dampco \in \Linfty$. 
		Then the following are equivalent:
		\begin{eqenumerate}
			\item \label{item:log stable}
			$\{ \semigroup[\gen[\dampco]] \}_{t \geq 0}$ is $1 / p$-logarithmically stable, that is, there exists $M > 0$ such that
			\begin{equation}
				\norm{ \semigroup[\gen[\dampco]] \gen[\dampco]^{-1} }_{ \sobolev{1} \times \hilbert \to \sobolev{1} \times \hilbert } 
				\leq M ( \log{(e + t)} )^{- 1 / p}
			\end{equation}
			holds for any $t \geq 0$.
			\item \label{item:exp resolvent of A}
			There exists $C > 0$ such that
			\begin{equation} 
				\norm{F}_{\sobolev{1} \times \hilbert} 
				\leq C e^{C\abs{\lambda}^{p}} \norm{ ( \gen[\dampco] - i \lambda ) F }_{\sobolev{1} \times \hilbert}
			\end{equation}
			holds for any $\lambda \in \R$ and $F \in \sobolev{2} \times \sobolev{1}$. 
		\end{eqenumerate}
	\end{quotetheorem}
	\begin{quotetheorem}[\cite{MR2460938}*{Theorem 1.1}] \label{theorem:o(1) stable = resolvent of A}
		Let $s > 0$ and $0 \leq \dampco \in \Linfty$. 
		Then the following are equivalent:
		\begin{eqenumerate}
			\item 	\label{item:o(1) stable}
			$\{ \semigroup[\gen[\dampco]] \}_{t \geq 0}$ is $o(1)$ stable, that is, 
			\begin{equation}
				\lim_{t \to \infty} \norm{ \semigroup[\gen[\dampco]] \gen[\dampco]^{-1} }_{ \sobolev{1} \times \hilbert \to \sobolev{1} \times \hilbert } 
				= 0
			\end{equation}
			holds.
			\item 
			For each fixed $\lambda \in \R$, there exists $C > 0$ such that
			\begin{equation} 
				\norm{F}_{\sobolev{1} \times \hilbert} \leq 
				C \norm{ ( \gen[\dampco] - i \lambda ) F }_{\sobolev{1} \times \hilbert}
			\end{equation}
			holds for any $F \in \sobolev{2} \times \sobolev{1}$. 
		\end{eqenumerate}
	\end{quotetheorem}
	Recent studies \cites{MR3547102, MR4119867, MR4143391, MR4363753, IS2022, MR3880304, MR4063985, MR3685285} gave several sufficient conditions and necessary conditions for the stabilities 
	using Theorem \ref{theorem:exp stable = uniform resolvent of A}, \ref{theorem:poly stable = poly resolvent of A}, \ref{theorem:log stable = exp resolvent of A}, \ref{theorem:o(1) stable = resolvent of A}.
	In order to state their results precisely, we introduce the so-called geometric control conditions (GCC).
	\begin{definition}
		Let $\Omega \subset \R^d$.
		\begin{description}
			\item[$d$-GCC] 
			$\Omega$ satisfies $d$-GCC when the following condition holds:
			there exists $r > 0$ such that 
			\begin{equation}
				\inf_{a \in \R^d} \measure{ B(a, r) \cap \Omega } > 0 , 
			\end{equation}
			where $B(a, r) \coloneqq \set{ x \in \R^d }{ \abs{x-a} < r }$.
			\item[$1$-GCC] 
			$\Omega$ satisfies $1$-GCC when the following condition holds:
			there exists $\ell > 0$ such that 
			\begin{equation}
				\inf_{\substack{ a \in \R^d ,  e \in S^{d-1}}} \mathcal{H}^{1}( L(a, e, \ell) \cap \Omega ) > 0 , 
			\end{equation}
			where $L(a, e, \ell) \coloneqq \set{a + \theta L e}{ 0 \leq \theta \leq 1 }$
			and $\mathcal{H}^{1}$ denotes the $1$-dimensional Hausdorff measure.
			\item[$0$-GCC] 
			$\Omega$ satisfies $0$-GCC when $\measure{\R^d \setminus \Omega} = 0$.
		\end{description}
	\end{definition}
	Actually we can also define $k$-GCC for any integers $2 \leq k \leq d-1$, 
	but we omit them here (see \cite{MR4363753}). 
	The papers we cited above show that 
	these geometric control conditions are closely related to the stabilities. 
	Here we give some of their results:
	\begin{quotetheorem}[\cite{MR4143391}*{Theorem 1}]
		Let $d = 1$, $s \geq 2$ and $0 \leq \dampco \in \Linfty$. 
		Then the following are equivalent:
		\begin{eqenumerate}
			\item[\eqref{item:exp stable}]  $\{ \semigroup[\gen[\dampco]] \}_{t \geq 0}$ is exponentially stable.
			\item There exists $\levelofdampco > 0$ such that $\set{x \in \R^d}{ \dampco(x) \geq \levelofdampco }$ satisfies $1$-GCC.
		\end{eqenumerate}
	\end{quotetheorem}
	\begin{quotetheorem}[\cite{MR4363753}*{Theorem 4}] \label{theorem:MR4363753 in intro}
		Let $d \geq 2$, $s \geq 2$ and $0 \leq \dampco \in \Linfty$ be uniformly continuous.
		Then \eqref{item:1-GCC} is sufficient for \eqref{item:exp stable}:
		\begin{eqenumerate}
			\item[\eqref{item:exp stable}] 
			$\{ \semigroup[\gen[\dampco]] \}_{t \geq 0}$ is exponentially stable.
			\item \label{item:1-GCC}
			There exists $\levelofdampco > 0$ such that $\set{x \in \R^d}{ \dampco(x) \geq \levelofdampco }$ satisfies $1$-GCC.
		\end{eqenumerate}
		Moreover, in the case $s = 2$, \eqref{item:exp stable} and \eqref{item:1-GCC} are equivalent.
	\end{quotetheorem}
	\begin{quotetheorem}[\cite{IS2022}*{Theorem 4}]
		Let $d \geq 1$, $0 < s < 2$ and $0 \leq \dampco \in \Linfty$.
		Then the following are equivalent:
		\begin{eqenumerate}
			\item[\eqref{item:exp stable}]
			$\{ \semigroup[\gen[\dampco]] \}_{t \geq 0}$ is exponentially stable.
			\item \label{item:0-GCC}
			There exists $\levelofdampco > 0$ such that $\set{x \in \R^d}{ \dampco(x) \geq \levelofdampco }$ satisfies $0$-GCC.
			In other words, $\essinf_{x \in \R^d} \dampco(x) > 0$ holds.
		\end{eqenumerate}
	\end{quotetheorem}
	\begin{quotetheorem}[\cite{IS2022}*{Theorem 5}]
		Let $d \geq 1$, $s > 0$ and $0 \leq \dampco \in \Linfty$.
		Then the following are equivalent:
		\begin{eqenumerate}
			\item[\eqref{item:o(1) stable}] 
			$\{ \semigroup[\gen[\dampco]] \}_{t \geq 0}$ is $o(1)$ stable.
			\item \label{item:d-GCC}
			There exists $\levelofdampco > 0$ such that $\set{x \in \R^d}{ \dampco(x) \geq \levelofdampco }$ satisfies $d$-GCC.
		\end{eqenumerate}
	\end{quotetheorem}
	See Table \ref{table:d=1, s < 2}, \ref{table:d=1, s geq 2}, \ref{table:d geq 2, s < 2}, \ref{table:d geq 2, s = 2}, \ref{table:d geq 2, s > 2}, \ref{table:d geq 2, s geq 4} and Section \ref{section:alternate proofs of some known results} for further results including the polynomial and logarithmic stabilities.
	These results suggest that ``the greater $s$ implies the more stability of $\{ \semigroup[\gen[\dampco]] \}_{t \geq 0}$''.
	We will prove this slogan precisely in the following sense, which is an extrapolation result respect to $s$:
	\begin{theorem} \label{theorem:stable for some s then for any s}
		Let $0 \leq \dampco \in \Linfty$, $s, s_0 > 0$, $p \geq 0$ and write 
		\begin{equation}
			\orderq \coloneqq (1 + p) s_0 / s - 1 .
		\end{equation}
		Then the following hold:
		\begin{eqenumerate}
			\item \label{item:exp stable for some s then}
			If $\{ \semigroup[\gen[\dampco][s_0]] \}_{t \geq 0}$ is exponentially stable, 
			then $\{ \semigroup[\gen[\dampco]] \}_{t \geq 0}$ is 	
			$( 2 \orderq[s][s_0][0] ) )^{-1}$-polynomially stable for any $0 < s < s_0$
			and exponentially stable for any $s \geq s_0$.
			\item \label{item:poly stable for some s then}
			If $\{ \semigroup[\gen[\dampco][s_0]] \}_{t \geq 0}$ is $1 / p$-polynomially stable, 
			then $\{ \semigroup[\gen[\dampco]] \}_{t \geq 0}$ is 
			$( 2 \orderq[s][s_0][ 2 p ] )^{-1}$-polynomially stable for any $0 < s < s_0$
			and $ s / ( 2 s_0 p ) $-polynomially stable for any $s \geq s_0$.
			\item \label{item:log stable for some s then}
			If $\{ \semigroup[\gen[\dampco][s_0]] \}_{t \geq 0}$ is $1 / p$-logarithmically stable, 
			then $\{ \semigroup[\gen[\dampco]] \}_{t \geq 0}$ is $s / ( s_0 p )$-logarithmically stable for any $s > 0$.
			\item \label{item:o(1) stable for some s then}
			If $\{ \semigroup[\gen[\dampco][s_0]] \}_{t \geq 0}$ is $o(1)$ stable, 
			then $\{ \semigroup[\gen[\dampco]] \}_{t \geq 0}$ is $o(1)$ stable for any $s > 0$.
		\end{eqenumerate}
	\end{theorem}
	In particular, Theorem \ref{theorem:stable for some s then for any s} allows us to obtain the polynomial stability for the non-fractional case $s = 2$
	via the exponential stability for sufficiently large $s > 2$.
	Furthermore, if we assume a reasonable assumption slightly stronger than the $1/p$-polynomial stability of $\{ \semigroup[\gen[\dampco][s_0]] \}_{t \geq 0}$, 
	then we can obtain the exponential stability of $\{ \semigroup[\gen[\dampco]] \}_{t \geq 0}$ for any $s \geq (1 + p) s_0$. 
	See Theorem \ref{theorem:stronger poly resolvent for some s0 implies exp stable for large s} for details. 
	
	Here recall that
	Theorem \ref{theorem:MR4363753 in intro} states that
	the condition \eqref{item:1-GCC}, 
	the existence of $\levelofdampco > 0$ such that $\set{x \in \R^d}{ \dampco(x) \geq \levelofdampco }$ satisfies $1$-GCC,
	is sufficient for the exponential stability of $\{ \semigroup[\gen[\dampco]] \}_{t \geq 0}$ in the case $s \geq 2$, 
	and it is also necessary in the case $s = 2$. 
	Since the greater $s$ implies the more stability, it is natural to ask:
	\begin{problem} \label{problem:large s}
		Let $d \geq 2$ and suppose that $s \geq 2$ is sufficiently large.
		Does there exists $0 \leq \dampco \in \Linfty \cap C^0(\R^d)$ satisfying both of the following conditions?
		\begin{eqenumerate}
			\item[\eqref{item:exp stable}] 
			$\{ \semigroup[\gen[\dampco]] \}_{t \geq 0}$ is exponentially stable.
			\item 
			For any $\levelofdampco > 0$, $\set{x \in \R^d}{ \dampco(x) \geq \levelofdampco }$ does not satisfy $1$-GCC.
		\end{eqenumerate}
	\end{problem}
	Though \cite{MR4119867} expected the existence, there are no such results as far as we know.
	We will show that actually it exists when $s \geq 4$ (see Theorem \ref{theorem:set S has a finite measure}, \ref{theorem:set S is contained in non-full periodic set}).
	
	Now we introduce an annihilating pair, which plays an important role in \cites{MR4143391, MR4363753, IS2022}.
	\begin{definition}
		Let $S, \Sigma \subset \R^d$.
		\begin{description}
			\item[Weak annihilation] 
			A pair $(S, \Sigma)$ is weakly annihilating when the following condition is satisfied:
			if $f \in \hilbert$ satisfies $\supp f \subset S$ and $\supp \Fourier[f] \subset \Sigma$, then $f = 0$.
			\item[Strong annihilation] 
			A pair $(S, \Sigma)$ is strongly annihilating when the following condition is satisfied:
			there exists $C > 0$ such that 
			\begin{equation} \label{eq:strong annihilation}
				\norm{f}_{\hilbert} \leq C ( \norm{\Fourier[f]}_{L^2(\R^d \setminus \Sigma)} + \norm{f}_{L^2(\R^d \setminus S)} )
			\end{equation} 
			holds for any $f \in \hilbert$.
			$\anihiconst[S][\Sigma]$ denotes the sharp constant $C$ of the strong annihilation inequality \eqref{eq:strong annihilation}. 
		\end{description}
	\end{definition}
	For example, the classical and well-known uncertainty principle states that $(S, \Sigma)$ is weakly annihilating if $S$ and $\Sigma$ are bounded. 
	In fact, the following stronger result is known: 
	if $\measure{S}, \measure{\Sigma} < \infty$, 
	then $(S, \Sigma)$ is strongly annihilating (see \cites{MR0461025, MR780328, MR2371612, MR1246419}).
	Also, the following Logvinenko--Serada theorem reveals the relation between $d$-GCC and the strong annihilation.
	\begin{quotetheorem}[\cite{MR0477719}] \label{theorem:L-S}
		Let $S \subset \R^d$. 
		Then the following are equivalent:
		\begin{eqenumerate}
			\item \label{item:L-S d-GCC}
			$\R^d \setminus S$ satisfies $d$-GCC.
			\item \label{item:L-S SA for any bdd sets}
			$(S, \Sigma)$ is strongly annihilating for any bounded measurable sets $\Sigma \subset \R^d$.
			\item \label{item:L-S SA for some non-null set}
			There exists $\Sigma \subset \R^d$ such that $(S, \Sigma)$ is strongly annihilating and $\measure{\Sigma} > 0$.
		\end{eqenumerate}
	\end{quotetheorem}
	We will establish an equivalent characterization of the exponential stability in terms of the strong annihilation. 
	Hereinafter, we write 
	\begin{gather} 
		\setS \coloneqq \set{ x \in \R^d }{ \dampco(x) < \levelofdampco } , \\
		\setSigma \coloneqq \set{ \xi \in \R^d }{ \abs{ (\abs{\xi}^2 + 1)^{s/4} - \lambda } < \levelofm } 
	\end{gather}
	for $s, \levelofdampco , \levelofm  > 0$, $\lambda \in \R$ and $0 \leq \dampco \in \Linfty$. 
	\begin{theorem} \label{theorem:exp stable = uniform strong annihilation}
		Let $s > 0$ and $0 \leq \dampco \in \Linfty$. 
		Then the following are equivalent:
		\begin{eqenumerate}
			\item[\eqref{item:exp stable}] 
			$\{ \semigroup[\gen[\dampco]] \}_{t \geq 0}$ is exponentially stable.
			\item \label{item:uniform strong annihilation}
			There exist $\levelofdampco, \levelofm > 0$ such that $\{ ( \setS, \setSigma ) \}_{\lambda \geq 0}$ is uniformly strongly annihilating, 
			that is, 
			there exists $C > 0$ such that 
			\begin{equation}
				\norm{f}_{\hilbert} \leq C ( \norm{\Fourier[f]}_{L^2(\R^d \setminus \setSigma )} + \norm{f}_{L^2(\R^d \setminus \setS)} )
			\end{equation}
			holds for any $\lambda \geq 0$ and $f \in \hilbert$.
		\end{eqenumerate}
	Furthermore, if \eqref{item:uniform strong annihilation} holds with $(C, \levelofdampco, \levelofm) = (C_0, \levelofdampco_0, \levelofm_0)$, 
	then we have
	\begin{equation}
	\norm{ \semigroup[\gen[\dampco]] }_{ \sobolev{1} \times \hilbert \to \sobolev{1} \times \hilbert } 
	\leq e^{ - \omega_0 t + \pi / 2} 
\end{equation}
	for any $t \geq 0$, where
	\begin{equation}
	\omega_0^{-1} = 8 ( 1 + \levelofdampco_0^{-1} + \dampconorm ) ( 1 + \levelofm_0^{-1} )^2 ( 1 + C_0 ) .
	\end{equation}
	\end{theorem}
	We can also prove similar equivalent characterizations of the polynomial, logarithmic and $o(1)$ stabilities.
	See Theorem \ref{theorem:poly stable = poly strong annihilation}, \ref{theorem:log stable = exp strong annihilation}, \ref{theorem:o(1) stable = pointwise strong annihilation} for details.
	
	Using Theorem \ref{theorem:stable for some s then for any s} and \ref{theorem:exp stable = uniform strong annihilation}, we obtain the following Theorem \ref{theorem:set S has a finite measure} and \ref{theorem:set S is contained in non-full periodic set}, new sufficient conditions for the exponential stability when $s$ is sufficiently large.
	\begin{theorem} \label{theorem:set S has a finite measure}
		Let $0 \leq \dampco \in \Linfty$
		and suppose that
		there exists $\levelofdampco > 0$ such that $\setS$ has a finite measure.
		Then $\{ \semigroup[\gen[\dampco]] \}_{t \geq 0}$ is 
		exponentially stable for any $s \geq 2d$
		and $( 4d/s - 2 )^{-1}$-polynomially stable for any $0 < s < 2d$.
	\end{theorem}
	\begin{theorem} \label{theorem:set S is contained in non-full periodic set}
		Let $0 \leq \dampco \in \Linfty$
		and suppose that
		there exist a closed periodic set $S \subsetneq \R^d$ and $\levelofdampco > 0$ such that $\setS \subset S$.
		Then $\{ \semigroup[\gen[\dampco]] \}_{t \geq 0}$ is 
		exponentially stable for any $s \geq 4$
		and $( 8/s - 2 )^{-1}$-polynomially stable for any $0 < s < 4$.
		Here `periodic' means that 
		there exists $\ell > 0$ satisfying
		\begin{equation}
			\set{ x + \ell n }{ x \in S , \, n \in \Z^d } = S .
		\end{equation}
	\end{theorem}
	We remark that 
	the the $1/2$-polynomial stability of $\{ \semigroup[\gen[\dampco][2]] \}_{t \geq 0}$ 
	under the assumption of Theorem \ref{theorem:set S is contained in non-full periodic set} 
	is already known by \cite{MR3685285}*{Theorem 1}.
	
	Finally, notice that both assumptions of Theorem \ref{theorem:set S has a finite measure} and \ref{theorem:set S is contained in non-full periodic set} are sufficient for 
	that there exists $\levelofdampco > 0$ such that $\R^d \setminus \setS$ satisfies $d$-GCC, 
	which implies the $s/2$-logarithmic stability of $\{ \semigroup[\gen[\dampco]] \}_{t \geq 0}$ for any $s > 0$ (see Theorem \ref{theorem:d-GCC log o(1)}).
	Therefore, the following seems to be a reasonable conjecture, though it remains open as of this writing.
	\begin{problem} \label{problem:conjecture}
		Let $d \geq 2$, $0 \leq \dampco \in \Linfty$ 
		and suppose that there exists $\levelofdampco > 0$ such that $\R^d \setminus \setS$ satisfies $d$-GCC.
		\begin{itemize}
			\item 
			Is $\{ \semigroup[\gen[\dampco]] \}_{t \geq 0}$ exponentially stable for sufficiently large $s > 2$?
			In particular, what about the case $s = 4$ or $s = 2d$?
			\item 
			If it is not, 
			then
			at least can we obtain the polynomial stability of $\{ \semigroup[\gen[\dampco]] \}_{t \geq 0}$?
		\end{itemize}
	\end{problem}
	This paper is organized as follows.
	
	In Section \ref{section:equivalence of the resolvent estimates}, 
	we will prove Theorem \ref{theorem:resolvent of A = resolvent of Laplacian},
	which gives the equivalence of the resolvent estimate of $\gen[\dampco]$; 
	\begin{equation} 
		\norm{F}_{\sobolev{1} \times \hilbert} 
		\leq C_0(\lambda) \norm{ ( \gen[\dampco] - i \lambda ) F }_{\sobolev{1} \times \hilbert}  ,
	\end{equation}
	the resolvent estimate of $\fracLaplacian[1]$ with an additional term $\norm{\dampco f}_{\hilbert}$;
	\begin{equation} \label{eq:intro resolvent of Laplacian}
		\norm{f}_{\hilbert} \leq C_1(\lambda) ( \norm{ ( \fracLaplacian[1] - \lambda ) f}_{\hilbert} + \norm{\dampco f}_{\hilbert} ) ,
	\end{equation} 
	and the resolvent estimate of $f \mapsto \dampco f$ for $f$ satisfying $\supp \Fourier[f] \subset \setSigma[\lambda][s][C_3(\lambda)^{-1}]$;
	\begin{equation} \label{eq:intro resolvent of B}
		\supp \Fourier[f] \subset \setSigma[\lambda][s][C_3(\lambda)^{-1}]
		\implies
		\norm{f}_{\hilbert} \leq C_4(\lambda) \norm{\dampco f}_{\hilbert} .
	\end{equation} 
	Combining Theorem \ref{theorem:resolvent of A = resolvent of Laplacian} with each one of the 
	Theorem \ref{theorem:exp stable = uniform resolvent of A}, \ref{theorem:poly stable = poly resolvent of A}, \ref{theorem:log stable = exp resolvent of A}, \ref{theorem:o(1) stable = resolvent of A}, 
	we obtain necessary and sufficient conditions for the stabilities of $\{ \semigroup[\gen[\dampco]] \}_{t \geq 0}$ 
	in terms of resolvent estimates of $\fracLaplacian[1]$ and $f \mapsto \dampco f$.
	See Theorem \ref{theorem:exp stable = uniform resolvent of Laplacian}, \ref{theorem:poly stable = poly resolvent of Laplacian}, \ref{theorem:log stable = exp resolvent of Laplacian}, \ref{theorem:o(1) stable = pointwise resolvent of Laplacian} for details.
	Our argument is based on \cite{MR4143391}*{Section 1} and \cite{MR4363753}*{Section 4}, 
	but we will prove it in a more systematic way.
	In fact, it holds even if 
	$f \mapsto \dampco f$ is replaced by any non-negative bounded operators on $\hilbert$. 
	In Section \ref{section:equivalence of the resolvent estimates} and \ref{section:stable for some s0 then for any s}, we will discuss in this setting.

	In Section \ref{section:stable for some s0 then for any s}, we prove Theorem \ref{theorem:stable for some s then for any s} by cases $s \geq s_0$ and $0 < s < s_0$ separately.
	We use \eqref{eq:intro resolvent of Laplacian} for $s \geq s_0$ and \eqref{eq:intro resolvent of B} for $0 < s < s_0$.

	In Section \ref{section:strong annihilation and resolvent estimates}, 
	at first we prove Theorem \ref{theorem:exp stable = uniform strong annihilation}
	and its analogies for the polynomial, logarithmic, $o(1)$ stabilities.
	Finally, using our results, we will prove Theorem \ref{theorem:set S has a finite measure} and \ref{theorem:set S is contained in non-full periodic set}.

	In Section \ref{section:alternate proofs of some known results}, 
	we give alternative proofs of some known results using our Theorem \ref{theorem:stable for some s then for any s} and \ref{theorem:exp stable = uniform strong annihilation}.
	\begin{table}[b]
		\begin{minipage}{0.5\textwidth}
			\caption{$d = 1$, $s < 2$}
			\centering
			\begin{tabular}{rcl} \label{table:d=1, s < 2}
				damping & & stability \\
				\hline
				$0$-GCC & $\Leftrightarrow$ & exponential \\
				$\Downarrow \quad \not \Uparrow$ & & $\Downarrow \quad \not \Uparrow$ \\
				$1$-GCC & $ \Leftrightarrow $ & polynomial \\
				& \rotatebox{-45}{$\Leftrightarrow$}  & $\Updownarrow$ \\
				& & $o(1)$\\
				& &
			\end{tabular}
		\end{minipage}%
		\begin{minipage}{0.5\textwidth}
			\caption{$d = 1, s \geq 2$} 
			\centering
			\begin{tabular}{rcl} \label{table:d=1, s geq 2}
				damping & & stability \\
				\hline
				$0$-GCC & & \\
				$\Downarrow \quad \not \Uparrow$ & & \\
				$1$-GCC & $\Leftrightarrow$ & exponential \\
				& \rotatebox{-45}{$\Leftrightarrow$}  & $\Updownarrow$ \\
				& & $o(1)$\\
				& &
			\end{tabular}
		\end{minipage}
		\begin{minipage}{0.5\textwidth}
			\caption{$d \geq 2, s < 2$}
			\centering
			\begin{tabular}{rcl} \label{table:d geq 2, s < 2}
				damping & & stability \\
				\hline
				$0$-GCC & $\Leftrightarrow$ & exponential \\
				$\Downarrow \quad \not \Uparrow$ & & $\Downarrow \quad \not \Uparrow$ \\
				$1$-GCC & $ \Rightarrow $ & polynomial \\
				$\Downarrow \quad \not \Uparrow$ & & $\Downarrow$ \\
				$d$-GCC & $ \Leftrightarrow $ & logarithmic \\
				& \rotatebox{-45}{$\Leftrightarrow$}  & $\Updownarrow$ \\
				& & $o(1)$\\
				& &
			\end{tabular}
		\end{minipage}%
		\begin{minipage}{0.5\linewidth}
			\caption{$d \geq 2, s = 2$}
			\centering
			\begin{tabular}{rcl} \label{table:d geq 2, s = 2}
				damping & & stability \\
				\hline
				$0$-GCC & & \\
				$\Downarrow \quad \not \Uparrow$ & & \\
				$1$-GCC & $\Leftrightarrow$ & exponential \\
				$\Downarrow \quad \not \Uparrow$ & & $\Downarrow \quad \not \Uparrow$ \\
				$d$-GCC & $ \Leftrightarrow $ & logarithmic \\
				& \rotatebox{-45}{$\Leftrightarrow$}  & $\Updownarrow$ \\
				& & $o(1)$\\
				& &
			\end{tabular}
		\end{minipage}
		\begin{minipage}{0.5\textwidth}
			\caption{$d \! \geq \! 2$, $2 \! < \! s \! \leq \! 4$}
			\centering
			\begin{tabular}{rcl} \label{table:d geq 2, s > 2}
				damping & & stability \\
				\hline
				$0$-GCC & & \\
				$\Downarrow \quad \not \Uparrow$ & & \\
				$1$-GCC &\begin{tabular}{@{}c@{}}
					$\Rightarrow$ \\
					\phantom{$\not \Leftarrow$}
				\end{tabular}& exponential \\
				$\Downarrow \quad \not \Uparrow$ & & $\Downarrow$ \\
				$d$-GCC & $ \Leftrightarrow $ & logarithmic \\
				& \rotatebox{-45}{$\Leftrightarrow$}  & $\Updownarrow$ \\
				& & $o(1)$\\
				& &
			\end{tabular}
		\end{minipage}%
		\begin{minipage}{0.5\textwidth}
			\caption{$d \geq 2, s \geq 4$}
			\centering
			\begin{tabular}{rcl} \label{table:d geq 2, s geq 4}
				damping & & stability \\
				\hline
				$0$-GCC & & \\
				$\Downarrow \quad \not \Uparrow$ & & \\
				$1$-GCC & \begin{tabular}{@{}c@{}}
					$\Rightarrow$ \\
					$\not \Leftarrow$
				\end{tabular} & exponential \\
				$\Downarrow \quad \not \Uparrow$ & & $\Downarrow$ \\
				$d$-GCC & $ \Leftrightarrow $ & logarithmic \\
				& \rotatebox{-45}{$\Leftrightarrow$}  & $\Updownarrow$ \\
				& & $o(1)$\\
				& &
			\end{tabular}
		\end{minipage}
		\begin{remark}
			Precisely, ``$1$-GCC $\implies$ polynomial'' in Table \ref{table:d geq 2, s < 2} and ``$1$-GCC $\implies$ exponential'' in Table \ref{table:d geq 2, s = 2}, \ref{table:d geq 2, s > 2}, \ref{table:d geq 2, s geq 4} need the uniformly continuity of damping coefficients.
			On the other hand, ``$1$-GCC $\impliedby$ exponential'' in Table \ref{table:d geq 2, s = 2} holds for continuous damping coefficients.
			See Theorem \ref{theorem:MR4363753 in intro}, \ref{theorem:d>1 1-GCC exp poly} and \ref{theorem:d>1 1-GCC exp poly with uniform continuity}.
		\end{remark}
	\end{table}
	
	\section{Equivalences of the resolvent estimates}
	\label{section:equivalence of the resolvent estimates}
	Let $\dampop \colon \hilbert \to \hilbert$ be bounded and consider the stability of
	\begin{equation} \label{eq:DFKGeq abstract}
		\begin{cases}
			u_{tt} (x, t) + \adjoint{\dampop} \dampop u_t(x, t) + \fracLaplacian[2] u(x, t) = 0 , & (x, t) \in \R^d \times \R_{> 0} , \\
			u(x, 0) = u_0(x) , & x \in \R^d , \\
			u_t(x, 0) = v_0(x) , & x \in \R^d .
		\end{cases}
	\end{equation}
	A typical example of $B$ is $f \mapsto \sqrt{ \dampco } f$ with $0 \leq \dampco \in \Linfty$, of course.
	We define $\gen \colon \sobolev{2} \times \sobolev{1} \to \sobolev{1} \to \hilbert$ by
	\begin{equation}
		\gen \coloneqq
		\begin{pmatrix}
			0 & 1 \\
			- \fracLaplacian[2] & - \adjoint{ \dampop } \dampop
		\end{pmatrix} 
	\end{equation}
	and prove the following Theorem \ref{theorem:resolvent of A = resolvent of Laplacian}:
	\begin{theorem} \label{theorem:resolvent of A = resolvent of Laplacian}
		Let $s > 0$, $\lambda \in \R$ and $\dampop \colon \hilbert \to \hilbert$ be bounded. 
		Then the following are equivalent:
		\begin{eqenumerate}
			\item \label{item:resolvent of A}
			There exists $C > 0$ such that
			\begin{equation} 
				\norm{F}_{\sobolev{1} \times \hilbert} \leq C \norm{ ( \gen - i \lambda ) F }_{\sobolev{1} \times \hilbert}
			\end{equation}
			holds for any $F \in \sobolev{2} \times \sobolev{1}$. 
			\item \label{item:resolvent of Laplacian}
			There exist $C_1, C_2 > 0$ such that 
			\begin{equation}
				\norm{f}_{\hilbert} 
				\leq C_1 \norm{ ( \fracLaplacian[1] - \abs{ \lambda } ) f}_{\hilbert} + C_2 \norm{\dampop f}_{\hilbert}
			\end{equation}
			holds for any $f \in \sobolev{1}$. 
			\item \label{item:resolvent of B}
			There exist $C, \levelofm > 0$ such that 
			\begin{equation}
				\norm{f}_{\hilbert} 
				\leq C \norm{\dampop f}_{\hilbert} 
			\end{equation}
			holds for any $f \in \hilbert$ satisfying $\supp \Fourier[f] \subset \setSigma[\abs{\lambda}]$.
		\end{eqenumerate}
		Furthermore, we have the following:
		\begin{itemize}
			\item 
			If \eqref{item:resolvent of A} holds with $C = C_0$, then 
			\eqref{item:resolvent of Laplacian} holds with
			\begin{equation}
				C_1 = C_0, \quad
				C_2 = C_0 \norm{ \dampop }_{\hilbert \to \hilbert} / \sqrt{2} .
			\end{equation}
			\item 
			If \eqref{item:resolvent of Laplacian} holds with $(C_1, C_2) = (C_{1, 0}, C_{2, 0})$, 
			then \eqref{item:resolvent of A} holds with 
			\begin{equation}
				C = 2 ( \sqrt{2} C_{3, 0} + 2 ( C_{2, 0} + C_{3, 0} \dampopnorm )^2 ) , 
			\end{equation}
		where
				\begin{equation}
			C_{3, 0} \coloneqq \max{ \{ C_{1, 0} , (1 + \abs{\lambda})^{-1} ( 1  +  C_{2, 0} \dampopnorm ) \} } . 
		\end{equation}
			\item
			If \eqref{item:resolvent of Laplacian} holds with $(C_1, C_2) = (C_{1, 0}, C_{2, 0})$, then 
			\eqref{item:resolvent of B} holds with
			\begin{equation}
				C = 2 C_{2, 0} , \quad \levelofm = 1/(2C_{1, 0}) .
			\end{equation}
			\item 
			If \eqref{item:resolvent of B} holds with $(C, \levelofm) = (C_0, \levelofm_0)$, 
			then \eqref{item:resolvent of Laplacian} holds with 
			\begin{equation}
				C_1 = \levelofm_0^{-1} ( 1 + C_0 \dampopnorm ) , \quad 
				C_2 = C_0 .
			\end{equation}
		\end{itemize}
	\end{theorem}
	Note that combining Theorem \ref{theorem:resolvent of A = resolvent of Laplacian} with each one of the 
	Theorem \ref{theorem:exp stable = uniform resolvent of A}, \ref{theorem:poly stable = poly resolvent of A}, \ref{theorem:log stable = exp resolvent of A}, \ref{theorem:o(1) stable = resolvent of A}, 
	we obtain following necessary and sufficient conditions for the stabilities of $\{ \semigroup \}_{t \geq 0}$.
	\begin{theorem} \label{theorem:exp stable = uniform resolvent of Laplacian}
		Let $s > 0$ and $\dampop \colon \hilbert \to \hilbert$ be bounded.
		Then the following are equivalent:
		\begin{eqenumerate}
			\item \label{item:exp stable dampop ver}
			$\{ \semigroup \}_{t \geq 0}$ is exponentially stable.
			\item \label{item:uniform resolvent of Laplacian}
			There exists $C > 0$ such that 
			\begin{equation}
				\norm{f}_{\hilbert} \leq C ( \norm{ ( \fracLaplacian[1] - \lambda ) f }_{\hilbert} + \norm{\dampop f}_{\hilbert} )
			\end{equation}
			holds for any $\lambda \geq 0$ and $f \in \sobolev{1}$.
			\item \label{item:uniform resolvent of B}
			There exist $C , \levelofm > 0$ such that 
			\begin{equation}
				\norm{f}_{\hilbert} \leq C \norm{\dampop f}_{\hilbert} 
			\end{equation}
			holds for any $\lambda \geq 0$ and $f \in \hilbert$ satisfying $\supp \Fourier[f] \subset \setSigma$.
		\end{eqenumerate}
	\end{theorem}
	\begin{theorem} \label{theorem:poly stable = poly resolvent of Laplacian}
		Let $s > 0$ and $\dampop \colon \hilbert \to \hilbert$ be bounded.
		Then the following are equivalent:
		\begin{eqenumerate}
			\item \label{item:poly stable dampop ver}
			There exists $p > 0$ such that $\{ \semigroup \}_{t \geq 0}$ is $1/p$-polynomially stable.
			\item \label{item:poly resolvent of Laplacian}
			There exist $C , \levelofm > 0$ and $p_1, p_2 \geq 0$ satisfying $p_1 + p_2 > 0$ such that 
			\begin{equation}
				\norm{f}_{\hilbert} \leq C (1 + \lambda)^{p_1} ( (1 + \lambda)^{p_2} \norm{ ( \fracLaplacian[1] - \lambda ) f }_{\hilbert} + \norm{\dampop f}_{\hilbert} )
			\end{equation}
			holds for any $\lambda \geq 0$ and $f \in \sobolev{1}$.
			\item \label{item:poly resolvent of B}
			There exist $C , \levelofm > 0$ and $p_1, p_2 \geq 0$ satisfying $p_1 + p_2 > 0$ such that 
			\begin{equation}
				\norm{f}_{\hilbert} \leq C (1 + \lambda)^{p_1}  \norm{\dampop f}_{\hilbert} 
			\end{equation}
			holds for any $\lambda \geq 0$ and $f \in \hilbert$ satisfying $\supp \Fourier[f] \subset \setSigma[\lambda][s][\levelofm (1 + \lambda)^{ - p_2}]$.
		\end{eqenumerate}
		Furthermore, we have the following:
		\begin{itemize}
			\item
			If \eqref{item:poly stable dampop ver} holds with $p = p_0$, 
			then \eqref{item:poly resolvent of Laplacian} and \eqref{item:poly resolvent of B} hold with 
			$(p_1, p_2) = (p_0, 0)$ and $(p_1, p_2) = (p_0, p_0)$, respectively.
			\item
			If \eqref{item:poly resolvent of Laplacian} or \eqref{item:poly resolvent of B} hold with $(p_1, p_2) = (p_{1, 0} , p_{2, 0})$, 
			then \eqref{item:poly stable dampop ver} holds with $p = 2 ( p_{1, 0} + p_{2, 0} ) $.
		\end{itemize}
	\end{theorem}
	\begin{theorem} \label{theorem:log stable = exp resolvent of Laplacian}
		Let $s , p > 0$ and $\dampop \colon \hilbert \to \hilbert$ be bounded.
		Then the following are equivalent:
		\begin{eqenumerate}
			\item \label{item:log stable dampop ver}
			$\{ \semigroup \}_{t \geq 0}$ is $1 / p$-logarithmically stable.
			\item \label{item:exp resolvent of Laplacian}
			There exist $C, \levelofm > 0$ and $p_1, p_2 \geq 0$ satisfying $\max \{ p_1 ,  p_2 \} = p$ such that 
			\begin{equation}
				\norm{f}_{\hilbert} 
				\leq C \exp{ ( C \lambda^{p_1} ) }
				( 
				\exp{ ( C \lambda^{p_2} ) } \norm{ ( \fracLaplacian[1] - \lambda ) f }_{\hilbert} 
				+ \norm{\dampop f}_{\hilbert} 
				)
			\end{equation}
			holds for any $\lambda \geq 0$ and $f \in \sobolev{1}$.
			\item \label{item:exp resolvent of B}
			There exist $C > 0$ and $p_1, p_2 \geq 0$ satisfying $\max \{ p_1 ,  p_2 \} = p$ such that 
			\begin{equation}
				\norm{f}_{\hilbert} \leq C \exp{ ( C \lambda^{p_1} ) } \norm{\dampop f}_{\hilbert} 
			\end{equation}
			holds for any $\lambda \geq 0$ and $f \in \hilbert$ satisfying $\supp \Fourier[f] \subset \setSigma[\lambda][s][\levelofm \exp{ ( -C \lambda^{p_2} ) }]$.
		\end{eqenumerate}
	\end{theorem}
	\begin{theorem} \label{theorem:o(1) stable = pointwise resolvent of Laplacian}
		Let $s > 0$ and $\dampop \colon \hilbert \to \hilbert$ be bounded.
		Then the following are equivalent:
		\begin{eqenumerate}
			\item \label{item:o(1) stable dampop ver} 
			$\{ \semigroup \}_{t \geq 0}$ is $o(1)$ stable.
			\item \label{item:pointwise resolvent of Laplacian dampop ver}
			For each fixed $\lambda \geq 0$, 
			there exists $C > 0$ such that 
			\begin{equation}
				\norm{f}_{\hilbert} \leq C ( \norm{ ( \fracLaplacian[1] - \lambda ) f }_{\hilbert} + \norm{\dampop f}_{\hilbert} )
			\end{equation}
			holds for any $f \in \sobolev{1}$.
			\item \label{item:pointwise resolvent of B}
			For each fixed $\lambda \geq 0$, 
			there exist $C, \levelofm > 0$ such that 
			\begin{equation}
				\norm{f}_{\hilbert} \leq C \norm{\dampop f}_{\hilbert} 
			\end{equation}
			holds for any $f \in \hilbert$ satisfying $\supp \Fourier[f] \subset \setSigma$.
		\end{eqenumerate}
	\end{theorem}
	\subsection{Proof of Theorem \ref{theorem:resolvent of A = resolvent of Laplacian}: \texorpdfstring{$\bf (\ref{item:resolvent of A}) \iff (\ref{item:resolvent of Laplacian})$}{(\ref{item:resolvent of A}) iff (\ref{item:resolvent of Laplacian})}}
	In this subsection, we will prove $\eqref{item:resolvent of A} \iff \eqref{item:resolvent of Laplacian}$ of Theorem \ref{theorem:resolvent of A = resolvent of Laplacian}.
	At first we prove Proposition \ref{prop:resolvent estimate step 1} and \ref{prop:resolvent estimate step 2}, 
	which lead us to $\eqref{item:resolvent of Laplacian} \implies \eqref{item:resolvent of A}$.
	\begin{proposition} \label{prop:resolvent estimate step 1}
		Let $C_1 , C_2 , s > 0$, $\lambda \in \R$ and $\dampop \colon \hilbert \to \hilbert$ be bounded. 
		Suppose that 
		\begin{equation} 
			\norm{f}_{\hilbert} \leq C_1 \norm{( \fracLaplacian[1] - \abs{ \lambda } )f}_{\hilbert} + C_2 \norm{\dampop f}_{\hilbert}
		\end{equation}
		holds for any $f \in \sobolev{1}$.
		Then we have
		\begin{align} 
			&\norm{F}_{\sobolev{1} \times \hilbert} \\
			&\leq
		\sqrt{2} C_3 \norm{ ( A(s, 0) - i \lambda ) F }_{\sobolev{1} \times \hilbert}  
		  + \sqrt{2} C_2 \norm{\dampop f_2}_{\hilbert}
		\end{align}
		for any $F = (f_1, f_2) \in \sobolev{2} \times \sobolev{1}$, where
		\begin{equation}
			C_3 \coloneqq \max{ \{ C_1 , (1 + \abs{\lambda})^{-1} ( 1  +  C_2 \dampopnorm ) \} } . 
		\end{equation}
	\end{proposition}
	\begin{proofof}{Proposition \ref{prop:resolvent estimate step 1}}
		Define 
		$\diagonizer{j} \colon \sobolev{(j+1)s / 2} \times \sobolev{js / 2} \to \sobolev{js / 2} \times \sobolev{js / 2}$ ($j = 0, 1$) by 
		\begin{equation}
			\diagonizer{j} \coloneqq \frac{1}{\sqrt{2}}
			\begin{pmatrix}
				\fracLaplacian[1] & - i \sgn{(\lambda)} \\
				\fracLaplacian[1] & i \sgn{(\lambda)}
			\end{pmatrix},
		\end{equation}
		where
		\begin{equation}
			\sgn{(\lambda)} = 
			\begin{cases}
				1 , & \lambda \geq 0 ,\\
				-1 , & \lambda < 0 .
			\end{cases}
		\end{equation}
		Then $\diagonizer{0}$ is unitary and $\diagonizer{1}$ diagonalizes $ \gen[0] - i \lambda $, that is,   
		\begin{gather}
			\diagonizer{0}^{-1}
			= \adjoint{ \diagonizer{0} } 
			= 
			\frac{1}{\sqrt{2}}
			\begin{pmatrix}
				\fracLaplacian[-] & \fracLaplacian[-] \\
				i \sgn{(\lambda)} & -i \sgn{(\lambda)}
			\end{pmatrix} , 
			\label{eq:P unitary} \\
			\begin{aligned} \label{eq:P diagonize}
				&\diagonizer{1} ( \gen[0] - i \lambda ) \diagonizer{1}^{-1} \\
				&= 
				\begin{pmatrix}
					i \sgn{(\lambda)} ( \fracLaplacian[1] - \abs{ \lambda }  )  & 0  \\
					0 & - i \sgn{(\lambda)} ( \fracLaplacian[1] + \abs{ \lambda }  ) 
				\end{pmatrix} .
			\end{aligned}			
		\end{gather}
		Fix $F = (f_1, f_2) \in \sobolev{2} \times \sobolev{1}$ 
		and write $G = (g_1, g_2) = \diagonizer{1} F \in \sobolev{1} \times \sobolev{1}$.
		Then we have
		\begin{gather} 
			\begin{aligned} \label{eq:resolvent estimate step 1-1}
				\norm{F}_{\sobolev{1} \times \hilbert}^2
				\overset{\eqref{eq:P unitary}}&{=} \norm{\diagonizer{1} F}_{\hilbert \times \hilbert}^2 \\
				&= \norm{g_1}_{\hilbert}^2 + \norm{g_2}_{\hilbert}^2 ,
			\end{aligned} \\
			\begin{aligned} \label{eq:resolvent estimate step 1-2}
				&\quad 
				\norm{ ( \gen[0] - i \lambda ) F }_{\sobolev{1} \times \hilbert}^2 \\
				\overset{\eqref{eq:P unitary}}&{=} 
				\norm{ \diagonizer{1} ( \gen[0] - i \lambda ) \diagonizer{1}^{-1} G }_{\hilbert \times \hilbert}^2 \\
				\overset{\eqref{eq:P diagonize}}&{=} 
				\Norm{ 
					\begin{pmatrix}
						i \sgn{(\lambda)} ( \fracLaplacian[1] - \abs{ \lambda }  )  & 0  \\
						0 & - i \sgn{(\lambda)} ( \fracLaplacian[1] + \abs{ \lambda }  ) 
					\end{pmatrix}
					\begin{pmatrix}
						g_1 \\
						g_2
				\end{pmatrix} }_{\hilbert \times \hilbert}^2 \\
				&= \norm{ ( \fracLaplacian[1] - \abs{ \lambda } ) g_1 }_{\hilbert}^2 + \norm{( \fracLaplacian[1] + \abs{ \lambda } ) g_2 }_{\hilbert}^2 . 
			\end{aligned}
		\end{gather}
		In addition, it is obvious that 
		\begin{equation} \label{eq:resolvent estimate step 1-3}
			\norm{g_2}_{\hilbert} \leq (1 + \abs{\lambda})^{-1} \norm{( \fracLaplacian[1] + \abs{ \lambda } ) g_2}_{\hilbert} 
		\end{equation}
		holds.
		Now, by the assumption of Proposition \ref{prop:resolvent estimate step 1}, we have
		\begin{equation} 
			\begin{aligned} \label{eq:resolvent estimate step 1-4}
				&\quad \norm{g_1}_{\hilbert} \\
				&\leq C_1 \norm{(\fracLaplacian[1] -  \abs{ \lambda } ) g_1}_{\hilbert} + C_2 \norm{\dampop g_1}_{\hilbert}  \\
				&\leq C_1 \norm{(\fracLaplacian[1] -  \abs{ \lambda } ) g_1}_{\hilbert} + C_2 \norm{\dampop ( g_1 - g_2 )}_{\hilbert} + C_2 \norm{ \dampop g_2 }_{\hilbert} \\
				&\leq C_1 \norm{(\fracLaplacian[1] -  \abs{ \lambda } ) g_1}_{\hilbert} + \sqrt{2} C_2 \norm{\dampop f_2}_{\hilbert} + C_2 \dampopnorm \norm{ g_2 }_{\hilbert}  .
			\end{aligned}
		\end{equation}
		Therefore, we obtain
		\begin{align}
			&\norm{F}_{\sobolev{1} \times \hilbert} \\
			&\overset{\eqref{eq:resolvent estimate step 1-1}}{\leq} 
			\norm{ g_1 }_{\hilbert} + \norm{g_2}_{\hilbert} \\
			&\overset{\eqref{eq:resolvent estimate step 1-4}}{\leq} 
			C_1 \norm{ ( \fracLaplacian[1] - \abs{ \lambda } ) g_1 }_{\hilbert} \\
			&\quad \quad + ( 1 + C_2 \dampopnorm ) \norm{ g_2 }_{\hilbert} + \sqrt{2} C_2 \norm{\dampop f_2}_{\hilbert} \\
			&\overset{ \eqref{eq:resolvent estimate step 1-3} }{\leq} 
			C_1 \norm{ ( \fracLaplacian[1] - \abs{ \lambda } ) g_1 }_{\hilbert} \\
			&\quad \quad + (1 + \abs{\lambda})^{-1} ( 1  +  C_2 \dampopnorm ) \norm{(\fracLaplacian[1] + \abs{ \lambda } ) g_2}_{\hilbert} \\
			&\quad \quad + \sqrt{2} C_2 \norm{\dampop f_2}_{\hilbert} \\
			&\overset{ \eqref{eq:resolvent estimate step 1-2} }{\leq} \sqrt{2} C_3 \norm{ ( A(s, 0) - i \lambda ) F }_{\sobolev{1} \times \hilbert} 
			 + \sqrt{2} C_2 \norm{\dampop f_2}_{\hilbert} , 
		\end{align}
	where
	\begin{equation}
	C_3 \coloneqq \max{ \{ C_1 , (1 + \abs{\lambda})^{-1} ( 1  +  C_2 \dampopnorm ) \} } . \qedhere
	\end{equation}
	\end{proofof}
	\begin{proposition} \label{prop:resolvent estimate step 2}
		Let $C_1, C_2, s > 0$, $\lambda \in \R$ and $\dampop \colon \hilbert \to \hilbert$ be bounded. 
		Then the inequality
		\begin{align} 
			&\quad \sqrt{2} C_1 \norm{ ( \gen[0] - i \lambda ) F }_{\sobolev{1} \times \hilbert} 
			+ \sqrt{2} C_2 \norm{\dampop f_2}_{\hilbert} \\
			&\leq  
	( \sqrt{2} C_1 + \delta^{-1} ( C_1 \dampopnorm + C_2 ) ) \norm{ ( \gen - i \lambda ) F }_{\sobolev{1} \times \hilbert} \\
	&\quad	+\delta ( C_1 \dampopnorm + C_2 ) \norm{ F }_{\sobolev{1} \times \hilbert} .
		\end{align}
		holds for any $\delta > 0$ and $F = (f_1, f_2) \in \sobolev{2} \times \sobolev{1}$.
	\end{proposition}
	\begin{proofof}{Proposition \ref{prop:resolvent estimate step 2}}
		Fix $\delta > 0$ and $F = (f_1, f_2) \in \sobolev{2} \times \sobolev{1}$.
		Since
		\begin{equation}
			\gen[0]
			= 
			\gen
			+ \begin{pmatrix}
				0 & 0 \\
				0 & \adjoint{\dampop} \dampop
			\end{pmatrix} ,
		\end{equation}
		we have
		\begin{gather}
			\begin{aligned} \label{eq:resolvent estimate step 2-1}
				&\quad \norm{ ( \gen[0] - i \lambda ) F }_{\sobolev{1} \times \hilbert} \\
				&\leq  
				\norm{ ( \gen - i \lambda ) F }_{\sobolev{1} \times \hilbert}
				+\norm{ \adjoint{\dampop} \dampop f_2 }_{\hilbert} ,
			\end{aligned} \\
			\begin{aligned}	\label{eq:resolvent estimate step 2-2}
				&\quad \innerproduct{ ( \gen[0] - i \lambda ) F, F }_{\sobolev{1} \times \hilbert} \\
				&= \innerproduct{ ( \gen - i \lambda ) F , F}_{\sobolev{1} \times \hilbert} 
				+ \innerproduct{ \adjoint{\dampop} \dampop f_2, f_2}_{\hilbert} .
			\end{aligned} 
		\end{gather}
		Then, since $\gen[0] - i \lambda$ is skew-adjoint, we get
		\begin{equation}
			\norm{ \dampop f_2 }_{\hilbert}^2
			= - \Re {\innerproduct{ ( \gen - i \lambda ) F, F }_{\sobolev{1} \times \hilbert} } 
		\end{equation}
		by taking the real part of \eqref{eq:resolvent estimate step 2-2}.
		Hence, we obtain
		\begin{equation} \label{eq:resolvent estimate step 2-3}
			\begin{aligned}
				\sqrt{2} \norm{ \dampop f_2 }_{\hilbert}
				&\leq \norm{ \delta^{-1} ( \gen - i \lambda ) F - \delta F }_{\sobolev{1} \times \hilbert} \\
				&\leq \delta^{-1} \norm{ ( \gen - i \lambda ) F }_{\sobolev{1} \times \hilbert} 
				+ \delta \norm{ F }_{\sobolev{1} \times \hilbert} .
			\end{aligned}
		\end{equation}
		Therefore, we conclude that the following inequality holds:
		\begin{align}
			& \sqrt{2} C_1 \norm{ ( A(s, 0) - i \lambda ) F }_{\sobolev{1} \times \hilbert} + \sqrt{2} C_2 \norm{\dampop f_2}_{\hilbert} \\
			\overset{\eqref{eq:resolvent estimate step 2-1}}&{\leq} \!\!
			 \sqrt{2} C_1 \norm{ ( \gen  -  i \lambda ) F }_{\sobolev{1} \times \hilbert} \\
			&\quad + \sqrt{2} ( C_1 \dampopnorm + C_2 ) \norm{ \dampop f_2 }_{\hilbert} \\
			\overset{\eqref{eq:resolvent estimate step 2-3}}&{\leq} \!\!
			( \sqrt{2} C_1 + \delta^{-1} ( C_1 \dampopnorm + C_2 ) ) \norm{ ( \gen - i \lambda ) F }_{\sobolev{1} \times \hilbert} \\
			&\quad	+\delta ( C_1 \dampopnorm + C_2 ) \norm{ F }_{\sobolev{1} \times \hilbert} . \qedhere
		\end{align}
	\end{proofof}
	\begin{proofof}{$\eqref{item:resolvent of Laplacian} \implies \eqref{item:resolvent of A}$}
		Suppose that the inequality 
	\begin{equation}
	\norm{f}_{\hilbert} 
	\leq C_1 \norm{ ( \fracLaplacian[1] - \abs{ \lambda } ) f}_{\hilbert} + C_2 \norm{\dampop f}_{\hilbert}
\end{equation}
		holds for any $f \in \sobolev{1}$. 
		Then Proposition \ref{prop:resolvent estimate step 1} and \ref{prop:resolvent estimate step 2} imply that 
		the inequality
		\begin{align}
			&\quad \norm{F}_{\sobolev{1} \times \hilbert} \\
			&\leq 
		( \sqrt{2} C_3 + \delta^{-1} ( C_2 + C_3 \dampopnorm ) ) \norm{ ( \gen - i \lambda ) F }_{\sobolev{1} \times \hilbert} \\
		&\quad + \delta ( C_2 + C_3 \dampopnorm ) \norm{ F }_{\sobolev{1} \times \hilbert} 
		\end{align}
		holds for any $\delta > 0$ and $F \in \sobolev{2} \times \sobolev{1}$, where
		\begin{equation}
		C_3 \coloneqq \max{ \{ C_1 , (1 + \abs{\lambda})^{-1} ( 1  +  C_2 \dampopnorm ) \} } . 
		\end{equation}
		Therefore, letting
		\begin{equation}
		\delta = \frac{1}{2 ( C_2 + C_3 \dampopnorm ) } ,
		\end{equation}
		we conclude that 
		\begin{align}
			&\quad \norm{F}_{\sobolev{1} \times \hilbert} \\
			&\leq
			2 ( \sqrt{2} C_3 + 2 ( C_2 + C_3 \dampopnorm )^2 ) \norm{ ( \gen - i \lambda ) F }_{\sobolev{1} \times \hilbert}
		\end{align}
		holds for any $F \in \sobolev{2} \times \sobolev{1}$. 
	\end{proofof}
	Next we prove $\eqref{item:resolvent of A} \implies \eqref{item:resolvent of Laplacian}$.
	\begin{proofof}{$\eqref{item:resolvent of A} \implies \eqref{item:resolvent of Laplacian}$}
		By the assumption \eqref{item:resolvent of A}, we can take $C > 0$ such that  
		\begin{equation} 
			\norm{F}_{\sobolev{1} \times \hilbert} \leq C \norm{ ( \gen - i \lambda ) F }_{\sobolev{1} \times \hilbert}
		\end{equation}
		holds for any $F \in \sobolev{2} \times \sobolev{1}$. 
		Fix $g \in \sobolev{1}$ and let
		$G = ( g, 0 ) \in \sobolev{1} \times \sobolev{1}$, $F = \diagonizer{1}^{-1} G \in \sobolev{2} \times \sobolev{1}$.
		Then we have
		\begin{gather}
			\norm{F}_{\sobolev{1} \times \hilbert}
			\overset{\eqref{eq:resolvent estimate step 1-1}}{=} 
			\norm{G}_{\hilbert \times \hilbert}
			= \norm{g}_{\hilbert} , \\
			\begin{aligned}
				&\quad \norm{ ( \gen - i \lambda ) F }_{\sobolev{1} \times \hilbert} \\
				&\leq  
				\norm{ ( \gen[0] - i \lambda ) F }_{\sobolev{1} \times \hilbert} \\
				&\quad + \frac{1}{\sqrt{2}}
				\Norm{ 
					\begin{pmatrix}
						0 & 0 \\
						0 & \adjoint{\dampop} \dampop
					\end{pmatrix}
					\begin{pmatrix}
						\fracLaplacian[-] & \fracLaplacian[-] \\
						i \sgn{(\lambda)} & -i \sgn{(\lambda)}
					\end{pmatrix} 
					\begin{pmatrix}
						g \\
						0
					\end{pmatrix} 
				}_{\sobolev{1} \times \hilbert} \\
				\overset{\eqref{eq:resolvent estimate step 1-2}}&=
				\norm{ ( \fracLaplacian[1] - \abs{ \lambda } ) g }_{\hilbert} + \norm{ \adjoint{\dampop} \dampop g }_{\hilbert}  / \sqrt{2} .
			\end{aligned}
		\end{gather}
		Therefore, we conclude that
		\begin{equation}
			\norm{g}_{\hilbert} \leq C ( \norm{ ( \fracLaplacian[1] - \abs{ \lambda } ) g }_{\hilbert} + \norm{ \dampop }_{\hilbert \to \hilbert} \norm{ B g }_{\hilbert}  / \sqrt{2} )
		\end{equation}
		holds for any $g \in \sobolev{1}$.
	\end{proofof}
	\subsection{Proof of Theorem \ref{theorem:resolvent of A = resolvent of Laplacian}: \texorpdfstring{$\bf (\ref{item:resolvent of Laplacian}) \iff (\ref{item:resolvent of B})$}{(\ref{item:resolvent of Laplacian}) iff (\ref{item:resolvent of B})}}
	In this subsection, we will prove $\eqref{item:resolvent of Laplacian} \iff \eqref{item:resolvent of B}$ of Theorem \ref{theorem:resolvent of A = resolvent of Laplacian}.
	We begin with the following Proposition \ref{prop:(B, Sigma) strong annihilation}.
	\begin{proposition} \label{prop:(B, Sigma) strong annihilation}
		Let $\Sigma \subset \R^d$ be measurable 
		and $\dampop \colon \hilbert \to \hilbert$ be bounded. 
		Then the following are equivalent:
		\begin{eqenumerate}
			\item 
			There exist $C_1, C_2 > 0$ such that 
			\begin{equation}
				\norm{f}_{\hilbert} \leq C_1 \norm{\Fourier[f]}_{L^2(\R^d \setminus \Sigma)} + C_2 \norm{\dampop f}_{\hilbert}
			\end{equation} holds for any $f \in \hilbert$.
			\label{item:(B, Sigma) strong annihilation 1}
			\item There exists $C > 0$ such that 
			\begin{equation}
				\norm{f}_{\hilbert} \leq C \norm{\dampop f}_{\hilbert}
			\end{equation}
			holds for any $f \in \hilbert$ satisfying $\supp \Fourier[f] \subset \Sigma$.
			\label{item:(B, Sigma) strong annihilation 2}
		\end{eqenumerate}
		Furthermore, we have the following:
		\begin{itemize}
			\item If \eqref{item:(B, Sigma) strong annihilation 1} holds with $(C_1, C_2) = (C_{1, 0}, C_{2, 0})$, 
			then \eqref{item:(B, Sigma) strong annihilation 2} holds with $C = C_{2, 0}$.
			\item If \eqref{item:(B, Sigma) strong annihilation 2} holds with $C = C_0$, 
			then \eqref{item:(B, Sigma) strong annihilation 1} holds with
			\begin{equation}
				C_1 = 1 + C_0 \dampopnorm , \quad C_2 = C_0 .
			\end{equation}
		\end{itemize}
	\end{proposition}
	\begin{proofof}{Proposition \ref{prop:(B, Sigma) strong annihilation}}
		Since $\eqref{item:(B, Sigma) strong annihilation 1} \implies \eqref{item:(B, Sigma) strong annihilation 2}$ is trivial, 
		it suffices to show $\eqref{item:(B, Sigma) strong annihilation 2} \implies \eqref{item:(B, Sigma) strong annihilation 1}$.
		Suppose that \eqref{item:(B, Sigma) strong annihilation 2} holds and
		define $P_\Sigma \colon \hilbert \to \hilbert$ by
		\begin{equation}
			\Fourier{P_\Sigma f} = \indicator{\Sigma} \Fourier[f] .
		\end{equation}
		Then we have
		\begin{align}
			\norm{f}_{\hilbert} 
			&\leq \norm{(I - P_\Sigma) f}_{\hilbert} 
			+ \norm{P_\Sigma f}_{\hilbert} \\
			\overset{\eqref{item:(B, Sigma) strong annihilation 2}}&\leq \norm{(I - P_\Sigma) f}_{\hilbert} 
			+ C \norm{\dampop P_\Sigma f}_{\hilbert} \\
			&\leq \norm{(I - P_\Sigma) f}_{\hilbert} 
			+ C ( \norm{\dampop (P_\Sigma - I) f}_{\hilbert} + \norm{\dampop f}_{\hilbert}  ) \\ 
			&\leq ( 1 + C \dampopnorm ) \norm{(I - P_\Sigma) f}_{\hilbert} + C \norm{\dampop f}_{\hilbert} \\
			&= ( 1 + C \dampopnorm ) \norm{\Fourier[f]}_{L^2(\R^d \setminus \Sigma)} + C \norm{\dampop f}_{\hilbert}
		\end{align}
		for any $f \in \hilbert$.
	\end{proofof}
	\begin{proposition} \label{prop:Fourier multiplier dampop ver}
		Let $\dampop \colon \hilbert \to \hilbert$ be bounded, $m \colon \R^d \to \C$ be measurable  
		and write
		\begin{gather}
			\setsigma \coloneqq \set{\xi \in \R^d}{ \abs{m(\xi)} < \levelofm } 
		\end{gather}
		for $\levelofm > 0$.
		Then the following are equivalent:
		\begin{eqenumerate}
			\item \label{item:Fourier multiplier dampop ver 2}
			There exist $C, \levelofm > 0$ such that 
			\begin{equation}
				\norm{f}_{\hilbert} \leq C \norm{\dampop f}_{\hilbert} 
			\end{equation}
			holds for any $f \in \hilbert$ satisfying $\supp \Fourier[f] \subset \setsigma$.
			\item \label{item:Fourier multiplier dampop ver 3}
			There exist $C_1, C_2 > 0$ such that 
			\begin{equation}
				\norm{f}_{\hilbert} \leq C_1 \norm{m \Fourier[f]}_{\hilbert} + C_2 \norm{\dampop f}_{\hilbert}
			\end{equation}
			holds for any $f \in \hilbert$.
		\end{eqenumerate}
		Furthermore, we have the following:
		\begin{itemize}
			\item 
			If \eqref{item:Fourier multiplier dampop ver 2} holds with $(C , \levelofm) = (C_0, \levelofm_0)$, then 
			\eqref{item:Fourier multiplier dampop ver 3} 
			holds with
			\begin{equation}
				C_1 = \levelofm_0^{-1} ( 1 + C_0 \dampopnorm ) , \quad
				C_2 = C_0
			\end{equation}
			\item
			If \eqref{item:Fourier multiplier dampop ver 3} holds with $(C_1, C_2) = (C_{1, 0} , C_{2, 0})$, then 
			\eqref{item:Fourier multiplier dampop ver 2}
			holds with
			\begin{equation}
				\quad C = 2 C_{2, 0} ,  \quad \levelofm = 1/(2C_{1, 0}) .
			\end{equation}
		\end{itemize}
	\end{proposition} 
	\begin{proofof}{Proposition \ref{prop:Fourier multiplier dampop ver}}
		$\eqref{item:Fourier multiplier dampop ver 2}\implies\eqref{item:Fourier multiplier dampop ver 3}$ is immediate from Proposition \ref{prop:(B, Sigma) strong annihilation} and
		\begin{gather}
			\norm{\Fourier[f]}_{L^2(\R^d \setminus \setsigma)} \leq \levelofm^{-1} \norm{m \Fourier[f]}_{\hilbert} .
		\end{gather} 
		Suppose that \eqref{item:Fourier multiplier dampop ver 3} holds
		and take $C_1, C_2 > 0$ such that 
		\begin{equation} 
			\norm{f}_{\hilbert} \leq C_1 \norm{ m \Fourier[ f ] }_{\hilbert} + C_2 \norm{\dampop f}_{\hilbert} 
		\end{equation}
		holds for any $f \in \hilbert$.
		Let $\levelofm = 1/(2C_1)$ and fix $f \in \hilbert$ satisfying $\supp \Fourier[f] \subset \setsigma$.
		Then we have
		\begin{equation}
			\norm{ m \Fourier[f] }_{\hilbert} \leq \levelofm \norm{f}_{\hilbert} = \norm{f}_{\hilbert} / (2 C_1) 
		\end{equation} 
		and thus
		\begin{equation} 
			\norm{f}_{\hilbert} \leq 
			2 ( C_1 \norm{ m \Fourier[f] }_{\hilbert} + C_2 \norm{\dampop f}_{\hilbert} ) - \norm{f}_{\hilbert} 
			\leq 2 C_2 \norm{\dampop f}_{\hilbert} 
		\end{equation}
		holds, as desired.
	\end{proofof}
	Applying Proposition \ref{prop:Fourier multiplier dampop ver} to
	\begin{equation}
		m_{\lambda, s}(\xi) = ( \abs{\xi}^2 + 1 )^{s/4} - \abs{ \lambda } ,
	\end{equation} 
	we obtain $\eqref{item:resolvent of Laplacian} \iff \eqref{item:resolvent of B}$ of Theorem \ref{theorem:resolvent of A = resolvent of Laplacian}.
	
	\section{Proof of Theorem \ref{theorem:stable for some s then for any s}}
	\label{section:stable for some s0 then for any s}
	In this section, We will prove Theorem \ref{theorem:stable for some s then for any s}. 
	Since proofs of \eqref{item:log stable for some s then}, \eqref{item:o(1) stable for some s then} are quite similar to those of \eqref{item:exp stable for some s then}, \eqref{item:poly stable for some s then}, 
	we only show \eqref{item:exp stable for some s then}, \eqref{item:poly stable for some s then} and omit details for \eqref{item:log stable for some s then}, \eqref{item:o(1) stable for some s then}.
	In addition, we also prove the following Theorem \ref{theorem:stronger poly resolvent for some s0 implies exp stable for large s}, 
	which requires a little stronger assumption and gives a much better result than \eqref{item:poly stable for some s then} of Theorem \ref{theorem:stable for some s then for any s}.
	\begin{theorem} \label{theorem:stronger poly resolvent for some s0 implies exp stable for large s}
		Let $C, s_0 , p > 0$ and $\dampop \colon \hilbert \to \hilbert$ be bounded.
		Suppose that the following inequality holds for any $\lambda \geq 0$ and $f \in \sobolev{0}$:
		\begin{equation}
			\norm{f}_{\hilbert} 
			\leq C
			\mleft(	
			( 1 + \lambda )^{ p } \norm{ ( \fracLaplacian[0] - \lambda ) f}_{\hilbert} 
			+ \norm{\dampop f}_{\hilbert} 
			\mright) .
		\end{equation}
		Then $\{ \semigroup \}_{t \geq 0}$ is
		$( 2 \orderq )^{-1}$-polynomially stable for any $0 < s < (1 + p) s_0$
		and exponentially stable for any $s \geq (1 + p) s_0$.
	\end{theorem}
	\subsection{In the case \texorpdfstring{$\bf s > s_0$}{s > s0}}
	\label{section:stable for some s0 then for any s > s0} 
	In this subsection, we will prove the following Proposition \ref{proposition:resolvent s0 to s > s0}:
	\begin{proposition} \label{proposition:resolvent s0 to s > s0}
		Let $C, s_0 > 0$, $p_1, p_2 \geq 0$ and $\dampop \colon \hilbert \to \hilbert$ be bounded.
		Suppose that the following inequality holds for any $\lambda \geq 0$ and $f \in \sobolev{0}$:
		\begin{equation}
			\norm{f}_{\hilbert} 
			\leq C
			\mleft(  \tfrac{ 1 + \lambda }{2} \mright)^{ p_1 } 
			\mleft(	
			\mleft(  \tfrac{ 1 + \lambda }{2} \mright)^{ p_2 } \norm{ ( \fracLaplacian[0] - \lambda ) f}_{\hilbert} 
			+ \norm{\dampop f}_{\hilbert} 
			\mright) .
		\end{equation}
		Then, for each fixed $s \geq s_0$, the inequality 
		\begin{equation}
			\norm{f}_{\hilbert} 
			\leq 	C \mleft( \tfrac{ 1 + \lambda }{2} \mright)^{ p_1 s_0 / s } 
			\mleft( 
			\mleft( \tfrac{ 1 + \lambda }{2} \mright)^{ \orderq[s][s_0][p_2] } \! \norm{ ( \fracLaplacian[1] - \lambda ) f }_{\hilbert}
			+ \norm{\dampop f}_{\hilbert}
			\mright) 
		\end{equation}
		holds for any $\lambda \geq 0$ and $f \in \sobolev{1}$.
	\end{proposition}
	Recall that 
	\begin{equation}
		\orderq = (1 + p) s_0 / s - 1 .
	\end{equation}
	Combining Proposition \ref{proposition:resolvent s0 to s > s0} with 
	Theorem \ref{theorem:exp stable = uniform resolvent of Laplacian}, \ref{theorem:poly stable = poly resolvent of Laplacian}, 
	we obtain \eqref{item:exp stable for some s then}, \eqref{item:poly stable for some s then} of Theorem \ref{theorem:stable for some s then for any s} 
	and Theorem \ref{theorem:stronger poly resolvent for some s0 implies exp stable for large s} in the case $s \geq s_0$.
	
	In order to prove Proposition \ref{proposition:resolvent s0 to s > s0}, 
	we need the following elementary Lemma \ref{lem:elementary inequality} and its consequence Corollary \ref{cor:cor of elementary inequality}.
	\begin{lemma} \label{lem:elementary inequality}
		Let $0 < r \leq 1$ and $a_1, a_2 \geq 0$ satisfy $(a_1, a_2) \ne (0, 0)$.
		Then we have
		\begin{equation} \label{eq:elementary inequality}
			\abs{ a_1^r - a_2^r } \leq 
			\mleft( \frac{ a + b }{2}  \mright)^{r - 1} \abs{ a_1 - a_2 } .
		\end{equation}
	\end{lemma}
	\begin{proofof}{Lemma \ref{lem:elementary inequality}}
		By the symmetricity, we assume that $a_1 \geq a_2 \geq 0$ holds without loss of generality.
		Furthermore, since the desired inequality is obvious if $a_1 > a_2 = 0$, 
		it is enough to prove the inequality in the case $a_1 \geq a_2 > 0$.
		In this case, we have
		\begin{equation}
			0 < a_1^{r - 1} \leq \mleft( \frac{a_1 + a_2}{2} \mright)^{r - 1} \leq a_2^{r - 1}
		\end{equation}
		and thus
		\begin{equation*}
			a_1^r - a_2^r = a_1^{r - 1} a_1 - a_2^{r - 1} a_2 \leq \mleft( \frac{a_1 + a_2}{2} \mright)^{r - 1} (a_1 - a_2) . \qedhere
		\end{equation*}
	\end{proofof}
	\begin{corollary} \label{cor:cor of elementary inequality}
		Let $s \geq s_0 > 0$. Then the following inequality
		\begin{equation} \label{eq:cor of elementary inequality 1}
			\norm{ ( \fracLaplacian[1][0] - \lambda^{ s_0 / s } ) f }_{\hilbert} 
			\leq 
			\mleft( \textstyle \frac{ 1 + \lambda }{2} \mright)^{s_0 / s - 1} \norm{ ( \fracLaplacian[1][0] - \lambda ) f }_{\hilbert} 
		\end{equation}
		holds for any $\lambda \geq 0$ and $f \in \sobolev{1}$.
	\end{corollary}
	\begin{proofof}{Corollary \ref{cor:cor of elementary inequality}}
		Fix $s \geq s_0 > 0$, $\lambda \geq 0$ and $\xi \in \R^d$. 
		Then
		using Lemma \ref{lem:elementary inequality} with $r = s_0 / s$, $a_1 = ( \abs{\xi}^2 + 1 )^{s / 4}$, $a_2 = \lambda$ implies that
		\begin{align}
			\abs{ ( \abs{\xi}^2 + 1 )^{s_0 / 4} - \lambda^{s_0 / s} } 
			\overset{\eqref{eq:elementary inequality}}&\leq \mleft( \frac{ ( \abs{\xi}^2 + 1 )^{s / 4} + \lambda }{2} \mright)^{s_0 / s - 1} \abs{ ( \abs{\xi}^2 + 1 )^{s / 4} - \lambda } \\
			&\leq \mleft( \frac{1 + \lambda }{2} \mright)^{s_0 / s - 1} \abs{ ( \abs{\xi}^2 + 1 )^{s / 4} - \lambda } 
		\end{align}
		holds.
		Now the desired inequality follows immediately.
	\end{proofof}
	\begin{proofof}{Proposition \ref{proposition:resolvent s0 to s > s0}}
		Suppose that the inequality
		\begin{equation}
			\norm{f}_{\hilbert} 
			\leq
			C 	\mleft( \tfrac{ 1 + \lambda_0 }{2} \mright)^{ p_1 }
			\mleft(	
			\mleft( \tfrac{ 1 + \lambda_0 }{2} \mright)^{ p_2 } \norm{ ( \fracLaplacian[0] - \lambda_0 ) f}_{\hilbert} 
			+ \norm{\dampop f}_{\hilbert} 
			\mright) 
		\end{equation}
		holds for any $\lambda_0 \geq 0$ and $f \in \sobolev{0}$.
		Fix $s \geq s_0$, $\lambda \geq 0$ and let $\lambda_0 = \lambda^{s_0/s}$. 
		Then, by Corollary \ref{cor:cor of elementary inequality}, we conclude that
		\begin{align}
			& \norm{f}_{\hilbert} \\
			&\leq
			C \mleft( \tfrac{ 1 + \lambda^{s_0 / s} }{2} \mright)^{ p_1 } 
			\mleft( 
			\mleft( \tfrac{ 1 + \lambda^{s_0/s} }{2} \mright)^{ p_2 } \norm{ ( \fracLaplacian[0] - \lambda^{s_0/s} ) f}_{\hilbert}
			+ \norm{\dampop f}_{\hilbert}
			\mright)  \\
			&\leq
			C \mleft(  \tfrac{ 1 + \lambda }{2} \mright)^{ p_1 s_0 / s } 
			\mleft( 
			\mleft(  \tfrac{ 1 + \lambda }{2} \mright)^{ \orderq[s][s_0][p_2] }  \norm{ ( \fracLaplacian[1] - \lambda ) f }_{\hilbert}
			+ \norm{\dampop f}_{\hilbert}
			\mright) 
		\end{align}
		holds for any $f \in \sobolev{1}$, as desired.
	\end{proofof}
	\subsection{In the case \texorpdfstring{$\bf s < s_0$}{s < s0}}
	\label{section:stable for some s0 then for any s < s0}
	In this subsection, we will prove the following Proposition \ref{prop:resolvent s0 to s < s0}.
	\begin{proposition} \label{prop:resolvent s0 to s < s0}
		Let $C , s_0 , \levelofm_0 > 0$, $p_1, p_2 \geq 0$ and $\dampop \colon \hilbert \to \hilbert$ be bounded. 
		Suppose that 
		the following inequality holds 
		for any $\lambda \geq 0$ and $f \in \hilbert$ satisfying $\supp \Fourier[f] \subset \setSigma[\lambda][s][ \levelofm_0 (1  + \lambda)^{- p_2} ]$:
		\begin{equation}
			\norm{f}_{ \hilbert } 
			\leq C
			( 1 + \lambda )^{ p_1 } \norm{\dampop f}_{ \hilbert } .
		\end{equation}
		Then, for each fixed $0 < s \leq s_0$, there exists $\levelofm > 0$ such that the inequality 
		\begin{equation}
			\norm{f}_{ \hilbert } 
			\leq C
			( 1 + \lambda )^{ p_1 s_0 / s }
			\norm{\dampop f}_{ \hilbert } 
		\end{equation}
		holds for any $\lambda \geq 0$ and $f \in \hilbert$ satisfying $\supp \Fourier[f] \subset \setSigma[\lambda][s][ \levelofm (1  + \lambda)^{- \orderq[s][s_0][p_2]} ]$.
	\end{proposition}
	Combining Proposition \ref{prop:resolvent s0 to s < s0} with 
	Theorem \ref{theorem:exp stable = uniform resolvent of Laplacian}, \ref{theorem:poly stable = poly resolvent of Laplacian}, 
	we obtain \eqref{item:exp stable for some s then}, \eqref{item:poly stable for some s then} of Theorem \ref{theorem:stable for some s then for any s} 
	and Theorem \ref{theorem:stronger poly resolvent for some s0 implies exp stable for large s} in the case $0 < s \leq s_0$.
	In order to prove Proposition \ref{proposition:resolvent s0 to s > s0}, 
we need the following elementary Lemma \ref{lem:elementary inequality 2} and its consequence Corollary \ref{cor:cor of elementary inequality 2}.
\begin{lemma} \label{lem:elementary inequality 2}
	Let $r > 1$ and $a_1, a_2 \geq 0$.
	Then we have
	\begin{equation} \label{eq:elementary inequality 2}
		\abs{ a_1^r - a_2^r } \leq 
		r \max\{a_1, a_2\}^{r - 1} \abs{ a_1 - a_2 } .
	\end{equation}
\end{lemma}
\begin{proofof}{Lemma \ref{lem:elementary inequality 2}}
This is immediate from the mean value theorem:
	\begin{equation}
	\abs{a_1^r - a_2^r} \leq \abs{a_1 - a_2} \sup_{0 \leq \theta \leq 1} r (\theta a_1 + (1 - \theta) a_2 )^{r - 1} = r \max\{a_1, a_2\}^{r - 1} \abs{ a_1 - a_2 } . \qedhere
	\end{equation}
\end{proofof}
	\begin{corollary} \label{cor:cor of elementary inequality 2}
	Let $0 < s < s_0$, $\levelofm_0 > 0$ and write
	\begin{equation}
	\levelofm \coloneqq \min\{ \levelofm_0 s / s_0, 1 \} .
	\end{equation}
	Then the inclusion
	\begin{equation} \label{eq:cor of elementary inequality 2}
	\setSigma[ \lambda ][s][ \levelofm (1 + \lambda)^{-\orderq} ] \subset \setSigma[ \lambda^{s_0 / s} ][s_0][ \levelofm_0 ( 1 + \lambda^{s_0 / s} )^{ - p } ]
	\end{equation}
	holds for any $p, \lambda \geq 0$.
\end{corollary}
\begin{proofof}{Corollary \ref{cor:cor of elementary inequality 2}}
	Fix $p, \lambda \geq 0$ and take $\xi \in \setSigma[ \lambda ][s][ \levelofm (1 + \lambda)^{-\orderq} ]$.
	Then, using Lemma \ref{lem:elementary inequality 2} with $r = s_0 / s$, $a_1 = ( \abs{\xi}^2 + 1 )^{s / 4}$, $a_2 = \lambda$ implies that
	\begin{align}
		\abs{ ( \abs{\xi}^2 + 1 )^{s_0 / 4} - \lambda^{s_0 / s} } 
		\overset{\eqref{eq:elementary inequality 2}}&\leq \frac{s_0}{s} \max\{ ( \abs{\xi}^2 + 1 )^{s / 4} , \lambda \}^{s_0 / s - 1} \abs{ ( \abs{\xi}^2 + 1 )^{s / 4} - \lambda } \\
		&\leq \frac{s_0}{s} ( \levelofm (1 + \lambda)^{-\orderq} + \lambda )^{s_0 / s - 1} \levelofm (1 + \lambda)^{-\orderq} \\
		&\leq \levelofm_0 ( 1 + \lambda )^{ - p s_0 / s} \\
		&\leq \levelofm_0 ( 1 + \lambda^{s_0 / s} )^{ - p } . \qedhere
	\end{align}
\end{proofof}
	\begin{proofof}{Proposition \ref{prop:resolvent s0 to s < s0}}
		Suppose that the inequality 
			\begin{equation}
			\norm{f}_{ \hilbert } 
			\leq C
			( 1 + \lambda_0 )^{ p_1 } 
			\norm{\dampop f}_{ \hilbert } 
		\end{equation}
		holds for any $\lambda_0 \geq 0$ and $f \in \hilbert$ satisfying $\supp \Fourier[f] \subset \setSigma[\lambda_0][s_0][ \levelofm_0 (1  + \lambda)^{- p_2} ]$.
		Fix $0 < s \leq s_0$, $\lambda \geq 0$, let 
		\begin{equation}
			\levelofm \coloneqq \min\{  \levelofm_0 s / s_0  , 1 \} , \quad
			\lambda_0 \coloneqq \lambda^{s_0 / s} ,
		\end{equation}
		and take $f \in \hilbert$ satisfying $\supp \Fourier[f] \subset \setSigma[\lambda][s][ \levelofm (1  + \lambda)^{- \orderq[s][s_0][p_2]} ]$.
		Then Corollary \ref{cor:cor of elementary inequality 2} implies $\supp \Fourier[f] \subset \setSigma[\lambda^{s_0 / s}][s_0][ \levelofm_0 (1  + \lambda^{s_0 / s})^{- p_2} ]$, 
		and thus we conclude that
		\begin{align*}
			\norm{f}_{ \hilbert } 
			&\leq C
			( 1 + \lambda^{s_0 / s} )^{ p_1 } 
			\norm{\dampop f}_{ \hilbert }   \\
			&\leq 
			C
			( 1 + \lambda )^{ p_1 s_0 / s }
			\norm{\dampop f}_{ \hilbert } . \qedhere
		\end{align*}
	\end{proofof}
	\section{The uncertainty principle and stabilities}
	\label{section:strong annihilation and resolvent estimates} 
	In this section, we give some necessary conditions and sufficient conditions for the stabilities. 
	At first we prove Theorem \ref{theorem:equivalence in the case 0 < s < 2 dampop ver}
	regarding the exponential stability in the case $0 < s < 2$.
	\begin{theorem} \label{theorem:equivalence in the case 0 < s < 2 dampop ver}
		Let $\dampop \colon \hilbert \to \hilbert$ be bounded and $\modulation{a}$ denotes $\modulation{a} f(x) = e^{i a \cdot x} f(x)$ for $a \in \R^d$.
		Then the following are equivalent:
		\begin{eqenumerate}
			\item \label{item:equivalence in the case 0 < s < 2 dampop ver 1}
			There exists $C > 0$ such that the inequality
			\begin{equation}
				\norm{f}_{\hilbert} \leq C \norm{\dampop f}_{\hilbert}
			\end{equation}
			for any $f \in \hilbert$.
			\item \label{item:equivalence in the case 0 < s < 2 dampop ver 2}
			$\{ \semigroup \}_{t \geq 0}$ is exponentially stable for any $s > 0$ 
			and there exists $C > 0$ such that the inequality
			\begin{equation}
				\norm{ \dampop f }_{\hilbert} \leq C \norm{ \dampop \modulation{a} f  }_{\hilbert}
			\end{equation}	
			holds for any $a \in \R^d$ and $f \in \hilbert$.
			\item \label{item:equivalence in the case 0 < s < 2 dampop ver 3}
			$\{ \semigroup \}_{t \geq 0}$ is exponentially stable for some $0 < s < 2$
			and there exists $C > 0$ such that the inequality
			\begin{equation}
				\norm{ \dampop f }_{\hilbert} \leq C \norm{ \dampop \modulation{a} f  }_{\hilbert}
			\end{equation}	
			holds for any $a \in \R^d$ and $f \in \hilbert$.
		\end{eqenumerate}	
	\end{theorem}
	Note that $\eqref{item:equivalence in the case 0 < s < 2 dampop ver 2}\implies\eqref{item:equivalence in the case 0 < s < 2 dampop ver 3}$ is trivial and
	$\eqref{item:equivalence in the case 0 < s < 2 dampop ver 1} \implies \eqref{item:equivalence in the case 0 < s < 2 dampop ver 2}$ is immediate 
	from $\eqref{item:uniform resolvent of B} \implies \eqref{item:exp stable dampop ver}$ of Theorem \ref{theorem:exp stable = uniform resolvent of Laplacian}
	and
	\begin{align}
		\norm{ \dampop f }_{\hilbert} 
		&\leq \dampopnorm \norm{ f }_{\hilbert} \\
		&= \dampopnorm \norm{ \modulation{a} f }_{\hilbert} \\
		\overset{\eqref{item:equivalence in the case 0 < s < 2 dampop ver 1}}
		&\leq C \dampopnorm \norm{\dampop \modulation{a} f}_{\hilbert}.
	\end{align}
	In order to prove $\eqref{item:equivalence in the case 0 < s < 2 dampop ver 3} \implies \eqref{item:equivalence in the case 0 < s < 2 dampop ver 1}$, 
	we need the following Proposition \ref{prop:setSigma 0 < s < 2}.
	\begin{proposition} \label{prop:setSigma 0 < s < 2}
		Let $0 < s < 2$, $\levelofm > 0$ and $\Omega \subset \R^d$ be bounded.
		Then there exist $a \in \R^d$ and $\lambda \geq 0$ such that $a + \Omega \subset \setSigma$ holds, 
		where $a + \Omega$ denotes
		\begin{equation}
			a + \Omega \coloneqq \set{ a + x }{ x \in \Omega } .
		\end{equation}
	\end{proposition}
	\begin{proofof}{Proposition \ref{prop:setSigma 0 < s < 2}}
		Without loss of generality, we assume that $\Omega$ is a ball centered at the origin.
		Fix $0 < s < 2$, $\levelofm , r > 0$ and let $\Omega = B(0, r)$, $\lambda \geq 1 + \levelofm$. 
		In this case, $\setSigma$ is an annulus as follows:
		\begin{align}
			\setSigma 
			&= \set{ \xi \in \R^d }{ \abs{ (\abs{\xi}^2 + 1)^{s/4} - \lambda } < \levelofm } \\
			&= \set{ \xi \in \R^d }{ ( (\lambda - \levelofm)^{4/s} - 1 )^{1/2} < \abs{\xi} < ( (\lambda + \levelofm)^{4/s} - 1 )^{1/2} } .
		\end{align}
		Since the mean value theorem and $0 < s < 2$ imply
		\begin{equation}
			\lim_{\lambda \to \infty} ( ( (\lambda + \levelofm)^{4/s} - 1 )^{1/2} - ( (\lambda - \levelofm)^{4/s} - 1 )^{1/2} ) 
			= \infty ,
		\end{equation}
		we can take $\lambda \geq 0$ such that
		\begin{equation}
			2 r \leq ( (\lambda + \levelofm)^{4/s} - 1 )^{1/2} - ( (\lambda - \levelofm)^{4/s} - 1 )^{1/2} .
		\end{equation}
		Therefore, taking $a \in \R^d$ satisfying
		\begin{equation}
			\abs{a} = \frac{ ( (\lambda + \levelofm)^{4/s} - 1 )^{1/2} + ( (\lambda - \levelofm)^{4/s} - 1 )^{1/2} }{2} ,
		\end{equation}
		we have $a + \Omega = B(a, r) \subset \setSigma$.
	\end{proofof}
	\begin{proofof}{$\eqref{item:equivalence in the case 0 < s < 2 dampop ver 3} \implies \eqref{item:equivalence in the case 0 < s < 2 dampop ver 1}$}
		Suppose that \eqref{item:equivalence in the case 0 < s < 2 dampop ver 3} holds.
		Then $\eqref{item:exp stable dampop ver} \implies \eqref{item:uniform resolvent of B}$ of Theorem \ref{theorem:exp stable = uniform resolvent of Laplacian} implies
		the existence of $C , \levelofm > 0$ and $0 < s < 2$ such that the inequality
		\begin{equation}
			\norm{f}_{\hilbert} \leq C \norm{\dampop \modulation{-a} f}_{\hilbert}
		\end{equation}
		holds for any $a \in \R^d$, $\lambda \geq 0$, and $f \in \hilbert$ satisfying $\supp \Fourier[f] \subset \setSigma$.
		
		Now fix $f \in \hilbert$ such that $\supp \Fourier[f]$ is compact.
		Then, by Proposition \ref{prop:setSigma 0 < s < 2}, 
		we can take $a \in \R^d$ and $\lambda \geq 0$ such that $\supp \Fourier[\modulation{a} f] = a + \supp \Fourier[f] \subset \setSigma$.
		Therefore, by the inequality above, 
		we conclude that
		\begin{equation*}
			\norm{f}_{\hilbert} = \norm{\modulation{a} f}_{\hilbert} \leq C \norm{\dampop \modulation{-a} \modulation{a} f}_{\hilbert} = C \norm{\dampop f}_{\hilbert} .
			\qedhere
		\end{equation*}
	\end{proofof}
	Hereinafter, we return to our original problem, that is, 
	the case of $\dampop \colon f \mapsto \sqrt{\dampco} f$ with $0 \leq \dampco \in \Linfty$. 
	Note that in this case Theorem \ref{theorem:equivalence in the case 0 < s < 2 dampop ver} is reduced to the following Corollary \ref{cor:equivalence in the case 0 < s < 2 dampco ver}.
	\begin{corollary}[\cite{IS2022}*{Theorem 4}] \label{cor:equivalence in the case 0 < s < 2 dampco ver}
		Let $0 \leq \dampco \in \Linfty$. 
		Then the following are equivalent:
		\begin{eqenumerate}
			\item[\eqref{item:0-GCC}]
			There exists $\levelofdampco > 0$ such that $\set{x \in \R^d}{ \dampco(x) \geq \levelofdampco }$ satisfies $0$-GCC.
			In other words, $\essinf_{x \in \R^d} \dampco(x) > 0$ holds.
			\item \label{item:equivalence in the case 0 < s < 2 dampco ver 2}
			$\{ \semigroup[\gen[\dampco]] \}_{t \geq 0}$ is exponentially stable for any $s > 0$.
			\item \label{item:equivalence in the case 0 < s < 2 dampco ver 3}
			$\{ \semigroup[\gen[\dampco]] \}_{t \geq 0}$ is exponentially stable for some $0 < s < 2$.
		\end{eqenumerate}	
	\end{corollary}
	Now we are going to prove Theorem \ref{theorem:exp stable = uniform strong annihilation}, 
	the equivalence of the uniform strong annihilation and the exponential stability.
	It is a consequence of Theorem \ref{theorem:exp stable = uniform resolvent of A} and the following Theorem \ref{theorem:resolvent of A = resolvent of Laplacian dampco ver}, which is a variant of Theorem \ref{theorem:resolvent of A = resolvent of Laplacian}.
		\begin{theorem} \label{theorem:resolvent of A = resolvent of Laplacian dampco ver}
		Let $s > 0$, $\lambda \in \R$ and $0 \leq \dampco \in \Linfty$. 
		Then the following are equivalent:
		\begin{eqenumerate}
			\item \label{item:resolvent of A dampco ver}
			There exists $C > 0$ such that
			\begin{equation} 
				\norm{F}_{\sobolev{1} \times \hilbert} \leq C \norm{ ( \gen[\dampco] - i \lambda ) F }_{\sobolev{1} \times \hilbert}
			\end{equation}
			holds for any $F \in \sobolev{2} \times \sobolev{1}$. 
			holds for any $f \in \sobolev{1}$. 
				\item \label{item:resolvent of Laplacian dampco ver}
			There exist $C_1, C_2, \levelofdampco > 0$ such that 
			\begin{equation}
				\norm{f}_{\hilbert} 
				\leq C_1 \norm{ ( \fracLaplacian[1] - \abs{ \lambda } ) f}_{\hilbert} + C_2 \norm{f}_{ L^2(\R^d \setminus \setS) }
			\end{equation}
		holds for any $f \in \sobolev{1}$.
					\item \label{item:resolvent of B dampco ver}
			There exist $\levelofdampco, \levelofm > 0$ such that 
			$( \setS, \setSigma[\abs{\lambda}] )$ is strongly annihilating.
		\end{eqenumerate}
		Furthermore, we have the following:
		\begin{itemize}
			\item 
			If \eqref{item:resolvent of Laplacian dampco ver} holds with $(C_1, C_2, \levelofdampco) = (C_{1, 0}, C_{2, 0}, \levelofdampco_0)$, 
			then \eqref{item:resolvent of A dampco ver} holds with 
			\begin{equation}
				C = 2 ( \sqrt{2} C_{3, 0} + 2 ( \levelofdampco^{-1/2} C_{2, 0} + C_{3, 0} \dampconorm^{1/2} )^2 ) , 
			\end{equation}
		where
			\begin{equation}
			C_{3, 0} \coloneqq \max{\{ C_{1, 0},  (1 + \abs{\lambda} )^{-1} ( 1 + C_{2, 0} ) \}} .
		\end{equation}
			\item
			If \eqref{item:resolvent of A dampco ver} holds with $C = C_0$, 
			then 
			\eqref{item:resolvent of B dampco ver} holds with
			\begin{equation}
				\anihiconst[\setS][\setSigma] \leq 2 ( 1 + \sqrt{2} C_0 \dampconorm ) , \quad 
				\levelofdampco = 1/(2 \sqrt{2} C), \quad 
				\levelofm = 1/(2 C) .
			\end{equation} 
			\item 
			If \eqref{item:resolvent of B dampco ver} holds with $(\levelofdampco, \levelofm) = (\levelofdampco_0, \levelofm_0)$, 
			then \eqref{item:resolvent of Laplacian dampco ver} and \eqref{item:resolvent of A dampco ver} hold with 
			\begin{equation}
				\levelofm_0 C_1 = C_2 = \anihiconst[\setS[\dampco][\levelofdampco_0]][\setSigma[\lambda][s][\levelofm_0]] , \quad 
				\levelofdampco = \levelofdampco_0 
			\end{equation}
		and
		\begin{equation}
		C = 8 (1 + \levelofm_0^{-1})^2 ( 1 + \levelofdampco_0^{-1} + \dampconorm ) (1 + \anihiconst[\setS[\dampco][\levelofdampco_0]][\setSigma[\lambda][s][\levelofm_0]] )^2 ,
		\end{equation}
	respectively.
		\end{itemize}
	\end{theorem}
	\begin{proofof}{Theorem \ref{theorem:resolvent of A = resolvent of Laplacian dampco ver}}
	Since $\eqref{item:resolvent of B dampco ver} \implies \eqref{item:resolvent of Laplacian dampco ver}$ is obvious, 
	it is enough to show
	$\eqref{item:resolvent of Laplacian dampco ver} \implies \eqref{item:resolvent of A dampco ver}$ 
	and $\eqref{item:resolvent of A dampco ver} \implies \eqref{item:resolvent of B dampco ver}$.
	
	Our proof of $\eqref{item:resolvent of Laplacian dampco ver} \implies \eqref{item:resolvent of A dampco ver}$ is basically the same as that of $\eqref{item:resolvent of Laplacian} \implies \eqref{item:resolvent of A}$ of Theorem \ref{theorem:resolvent of A = resolvent of Laplacian}.
	Suppose that 
	\begin{equation}
		\norm{f}_{\hilbert} 
		\leq C_1 \norm{ ( \fracLaplacian[1] - \abs{ \lambda } ) f}_{\hilbert} + C_2 \norm{f}_{ L^2(\R^d \setminus \setS) }
	\end{equation}
	holds for any $f \in \sobolev{1}$.
	Then using Proposition \ref{prop:resolvent estimate step 1} with $\dampop \colon f \mapsto \indicator{\R^d \setminus \setS} f$ implies that
	the inequality  
		\begin{align}
		&\quad \norm{F}_{\sobolev{1} \times \hilbert} \\
		&\leq \sqrt{2} C_3 \norm{ ( A(s, 0) - i \lambda ) F }_{\sobolev{1} \times \hilbert} + \sqrt{2} C_2 \norm{f}_{ L^2(\R^d \setminus \setS) } \\
		&\leq \sqrt{2} C_3 \norm{ ( A(s, 0) - i \lambda ) F }_{\sobolev{1} \times \hilbert} + \sqrt{2} \levelofdampco^{-1/2} C_2 \norm{\sqrt{\dampco} f}_{ \hilbert } 
	\end{align}
	holds for any $F \in \sobolev{2} \times \sobolev{1}$, where
	\begin{equation}
	C_3 \coloneqq \max{\{ C_{1},  (1 + \abs{\lambda} )^{-1} ( 1 + C_{2} ) \}} .
	\end{equation}
	Hence, using Proposition \ref{prop:resolvent estimate step 2} with $\dampop \colon f \mapsto \sqrt{\dampco} f$, we obtain the inequality 
		\begin{align}
		&\quad \norm{F}_{\sobolev{1} \times \hilbert} \\
		&\leq 
			( \sqrt{2} C_3 + \delta^{-1} ( \levelofdampco^{-1/2} C_2 + C_3 \dampconorm^{1/2} ) ) \norm{ ( \gen - i \lambda ) F }_{\sobolev{1} \times \hilbert} \\
	&\quad + \delta ( \levelofdampco^{-1/2} C_2 + C_3 \dampconorm^{1/2} ) \norm{ F }_{\sobolev{1} \times \hilbert} 
	\end{align}
	for any $\delta > 0$ and $F \in \sobolev{2} \times \sobolev{1}$.
	Therefore, letting
	\begin{equation}
		\delta = \frac{1}{2 ( \levelofdampco^{-1/2} C_2 + C_3 \dampconorm^{1/2} ) } ,
	\end{equation}
	we conclude that 
	\begin{align}
		&\quad \norm{F}_{\sobolev{1} \times \hilbert} \\
		&\leq
		2 ( \sqrt{2} C_3 + 2 ( \levelofdampco^{-1/2} C_2 + C_3 \dampconorm^{1/2} )^2 ) \norm{ ( \gen - i \lambda ) F }_{\sobolev{1} \times \hilbert} 
	\end{align}
	holds for any $F \in \sobolev{2} \times \sobolev{1}$. 
	
	Next we prove $\eqref{item:resolvent of A dampco ver} \implies \eqref{item:resolvent of B dampco ver}$. 
	Suppose that 
		\begin{equation} 
		\norm{F}_{\sobolev{1} \times \hilbert} \leq C \norm{ ( \gen[\dampco] - i \lambda ) F }_{\sobolev{1} \times \hilbert}
	\end{equation}
	holds for any $F \in \sobolev{2} \times \sobolev{1}$. 
	By the same argument as in the proof of $\eqref{item:resolvent of A} \implies \eqref{item:resolvent of B}$ of Theorem \ref{theorem:resolvent of A = resolvent of Laplacian}, we have the inequality
		\begin{equation}
		\norm{f}_{\hilbert} \leq C ( \norm{ ( \fracLaplacian[1] - \abs{ \lambda } ) f }_{\hilbert} + \norm{ \dampco f }_{\hilbert}  / \sqrt{2} )
	\end{equation}
	for any $f \in \sobolev{1}$.
	Then $\eqref{item:Fourier multiplier dampop ver 3} \implies \eqref{item:Fourier multiplier dampop ver 2}$ of Proposition \ref{prop:Fourier multiplier dampop ver} implies that 
	\begin{equation}
		\norm{f}_{\hilbert} 
		\leq \sqrt{2} C \norm{ \dampco f}_{\hilbert} 
	\end{equation}
	holds for any $f \in \hilbert$ satisfying $\supp \Fourier[f] \subset \setSigma$, where $\levelofm = 1/(2 C)$.
	
	Next we use $\eqref{item:(B, Sigma) strong annihilation 2} \implies \eqref{item:(B, Sigma) strong annihilation 1}$ of Proposition \ref{prop:(B, Sigma) strong annihilation} and obtain 
	\begin{equation}
		\norm{f}_{\hilbert} 
		\leq ( 1 + \sqrt{2} C \dampconorm ) \norm{ \Fourier[f] }_{L^2(\R^d \setminus \setSigma)} + \sqrt{2} C \norm{\dampco f}_{\hilbert} 
	\end{equation}
	for any $f \in \hilbert$.
	Now fix $f \in \hilbert$ and let $\levelofdampco = 1 / (2 \sqrt{2} C)$.
	Then we have
	\begin{align}
		\norm{ \dampco f }_{\hilbert} 
		&\leq \norm{ \dampco f }_{L^2(\setS)} + \norm{ \dampco f }_{L^2(\R^d \setminus \setS)} \\
		&\leq (2 \sqrt{2} C)^{-1} \norm{f}_{\hilbert} + \dampconorm \norm{ f }_{L^2(\R^d \setminus \setS)} 
	\end{align}
	and thus
	\begin{align}
	&\quad \norm{f}_{\hilbert} \\
	&\leq 2 ( 1 + \sqrt{2} C \dampconorm ) \norm{ \Fourier[f] }_{L^2(\R^d \setminus \setSigma)} + 2 \sqrt{2} C \dampconorm \norm{ f }_{L^2(\R^d \setminus \setS)} \\
	&\leq 2 ( 1 + \sqrt{2} C \dampconorm ) ( \norm{ \Fourier[f] }_{L^2(\R^d \setminus \setSigma)} + \norm{ f }_{L^2(\R^d \setminus \setS)}  ) . \qedhere
\end{align}
	\end{proofof}
	Combining Theorem \ref{theorem:exp stable = uniform resolvent of A} with Theorem \ref{theorem:poly stable = poly resolvent of Laplacian}, \ref{theorem:log stable = exp resolvent of Laplacian}, \ref{theorem:o(1) stable = pointwise resolvent of Laplacian}, 
	we also obtain similar results for the polynomial, logarithmic, $o(1)$ stabilities as follows. 
	\begin{theorem} \label{theorem:poly stable = poly strong annihilation}
		Let $s > 0$ and $0 \leq \dampco \in \Linfty$.
		Then the following are equivalent:
		\begin{eqenumerate}
			\item[\eqref{item:poly stable}]
			There exists $p > 0$ such that $\{ \semigroup[\gen[\dampco]] \}_{t \geq 0}$ is $1/p$-polynomially stable.
			\item \label{item:poly strong annihilation}
			There exist $C, \levelofdampco, \levelofm > 0$ and $p_1, p_2 , p_3 \geq 0$ satisfying $p_1 + p_2 + p_3 > 0$ such that 
			\begin{equation}
				( \setS[\dampco][ \levelofdampco (1 + \lambda)^{- p_1 } ] , \setSigma[\lambda][s][ \levelofm (1 + \lambda)^{ - p_2} ] )
			\end{equation}
			is strongly annihilating and
			\begin{equation}
				\anihiconst[ \setS[\dampco][ \levelofdampco (1 + \lambda)^{- p_1 } ] ][ \setSigma[\lambda][s][ \levelofm (1 + \lambda)^{ - p_2} ] ] \leq C (1 + \lambda)^{p_3} 
			\end{equation}
			holds for any $\lambda \geq 0$.
		\end{eqenumerate}
		Furthermore, we have the following:
		\begin{itemize}
			\item
			If \eqref{item:poly stable} holds with 
			$p = p_0$, 
			then \eqref{item:poly strong annihilation} holds with 
			\begin{equation}
				p_1 = p_2 = p_3 = p_0.
			\end{equation}
			\item
			If \eqref{item:poly strong annihilation} holds with 
			$(p_1, p_2, p_3) = (p_{1, 0} , p_{2, 0}, p_{3, 0})$, 
			then \eqref{item:poly stable} holds with 
			\begin{equation}
				p =  \max \{ p_{1, 0} , 2 p_{2, 0} \} + 2 p_{3, 0} .
			\end{equation}
		\end{itemize}
	\end{theorem}
	\begin{theorem} \label{theorem:log stable = exp strong annihilation}
		Let $s , p > 0$ and $0 \leq \dampco \in \Linfty$.
		Then the following are equivalent:
		\begin{eqenumerate}
			\item[\eqref{item:log stable}]
			$\{ \semigroup[\gen[\dampco]] \}_{t \geq 0}$ is $1 / p$-logarithmically stable.
			\item \label{item:exp strong annihilation}
			There exist $C, \levelofdampco, \levelofm > 0$ and $p_1, p_2 , p_3 \geq 0$ satisfying $\max \{ p_1 , p_2 , p_3 \} = p$ such that 
			\begin{equation}
				( \setS[\dampco][ \levelofdampco \exp{( - C \lambda^{p_2} )} ] , \setSigma[\lambda][s][ \levelofm \exp{( - C \lambda^{p_3} )} ] )
			\end{equation}
			is strongly annihilating and
			\begin{equation}
				\anihiconst[ \setS[\dampco][ \levelofdampco \exp{( - C \lambda^{ p_2} )} ] ][ \setSigma[\lambda][s][ \levelofm \exp{( - C \lambda^{ p_3} )} ] ] \leq C \exp{( C \lambda^{p_1} )}
			\end{equation}
			holds for any $\lambda \geq 0$.
		\end{eqenumerate}
	\end{theorem}
	\begin{theorem} \label{theorem:o(1) stable = pointwise strong annihilation}
		Let $s > 0$ and $0 \leq \dampco \in \Linfty$.
		Then the following are equivalent:
		\begin{eqenumerate}
			\item[\eqref{item:o(1) stable}]
			$\{ \semigroup[\gen[\dampco]] \}_{t \geq 0}$ is $o(1)$ stable.
			\item \label{item:pointwise strong annihilation}
			For any $\lambda \geq 0$, 
			there exist $\levelofdampco, \levelofm > 0$ such that $( \setS, \setSigma )$ is strongly annihilating.
		\end{eqenumerate}
	\end{theorem}
	Finally, 
	we prove Theorem \ref{theorem:set S has a finite measure} and \ref{theorem:set S is contained in non-full periodic set} using the following Theorem \ref{thm:uniform pair 1} and \ref{theorem:set S is contained in non-full periodic set s=4}, respectively:
	\begin{quotetheorem}[\cite{MR1246419} for $d=1$, \cite{MR2371612} for $d \geq 2$] \label{thm:uniform pair 1}
		If $S, \Sigma \subset \R^d$ satisfy $\measure{S}, \measure{\Sigma} < \infty$, 
		then 
		$(S, \Sigma)$ is strongly annihilating
		and $\anihiconst[S][\Sigma] \leq C_d \exp{( C_d \measure{S} \measure{\Sigma} )}$, 
		where $C_d > 0$ denotes a constant depending only on $d$.
	\end{quotetheorem}
	\begin{quotetheorem}[\cite{MR3685285}*{Theorem 2}] \label{theorem:set S is contained in non-full periodic set s=4}
		Let $S \subsetneq \R^d$ be a closed periodic set and $\levelofdampco > 0$.
		Then there exists $C > 0$ such that the inequality
		\begin{equation}
			\norm{f}_{\hilbert} 
			\leq C ( 
			\norm{ ( ( - \Laplacian + 1) - \lambda ) f }_{\hilbert} 
			+ \norm{ f }_{L^2( \R^d \setminus S )}
			)
		\end{equation}
		holds for any $\lambda \geq 0$ and $f \in \sobolev{{2}}$.
	\end{quotetheorem}
	\begin{proofof}{Theorem \ref{theorem:set S has a finite measure}}
		Take $\levelofdampco > 0$ such that $\measure{\setS}$ is finite.
		By Theorem \ref{theorem:stable for some s then for any s} and \ref{theorem:exp stable = uniform strong annihilation}, 
		it is enough to prove that $\{ \setS, \setSigma[\lambda][2d][1] \}_{\lambda \geq 0}$ is uniformly strongly annihilating.
		Moreover, by Theorem \ref{thm:uniform pair 1} and the continuity of $\lambda \mapsto \measure{\setSigma[\lambda][2d][1]}$, 
		it suffices to show that 
		\begin{equation}
			\limsup_{\lambda \to \infty} { \measure{\setSigma[\lambda][2d][1]} } < \infty .
		\end{equation}
		Let $\lambda \geq 2$. 
		Then we have
		\begin{align}
			\setSigma[\lambda][2d][1] 
			&= \set{ \xi \in \R^d }{ \abs{ ( \abs{\xi}^2 + 1 )^{d/2} - \lambda } < 1 } \\
			&= \set{ \xi \in \R^d }{ ( (\lambda - 1)^{2/d} - 1 )^{1/2} < \abs{\xi} < ( (\lambda + 1)^{2/d} - 1 )^{1/2} } 
		\end{align}
		and thus
		\begin{equation}
			\measure{ \setSigma[\lambda][2d][1] }
			= ( ( (\lambda + 1)^{2/d} - 1 )^{d/2} - ( (\lambda - 1)^{2/d} - 1 )^{d/2} ) \measure{B(0, 1)} .
		\end{equation}
		Therefore, using the mean value theorem, we conclude that the desired result holds.
	\end{proofof}
	\begin{proofof}{Theorem \ref{theorem:set S is contained in non-full periodic set}}
		By the assumption of Theorem \ref{theorem:set S is contained in non-full periodic set}, 
		there exist a closed periodic set $S \subsetneq \R^d$ and $\levelofdampco > 0$ such that $\setS \subset S$.
		Then, by Theorem \ref{theorem:set S is contained in non-full periodic set s=4}, 
		there exists $C > 0$ such that the inequality
		\begin{equation}
			\norm{f}_{\hilbert} 
			\leq C ( 
			\norm{ ( ( - \Laplacian + 1) - \lambda ) f }_{\hilbert} 
			+ \norm{ f }_{ L^2(\R^d \setminus \setS) }
			)
		\end{equation}
		holds for any $\lambda \geq 0$ and $f \in \sobolev{{2}}$.
		Hence, by Theorem \ref{theorem:exp stable = uniform resolvent of A}, \ref{theorem:stable for some s then for any s} and \ref{theorem:resolvent of A = resolvent of Laplacian dampco ver}, 
		we conclude that $\{ \semigroup[\gen[\dampco]] \}_{t \geq 0}$ is 
		exponentially stable for any $s \geq 4$ 
		and $s/(8 - 2s)$-polynomially stable for any $0 < s < 4$.
	\end{proofof}
	\section{Alternate proofs of some known results}
	\label{section:alternate proofs of some known results}
	In this section, we give alternate proofs of some known results.
	At first we prove Theorem \ref{theorem:d-GCC log o(1)}, which gives a necessary and sufficient condition for the logarithmic and $o(1)$ stabilities. 
	\begin{quotetheorem} \label{theorem:d-GCC log o(1)}
		Let $0 \leq \dampco \in \Linfty$.
		Then the following are equivalent:
		\begin{eqenumerate}
			\item[\eqref{item:d-GCC}]
			There exists $\levelofdampco > 0$ such that $\R^d \setminus \setS$ satisfies $d$-GCC.
			\item \label{item:log stable for any s > 0}
			For any $s > 0$, $\{ \semigroup[\gen[\dampco]] \}_{t \geq 0}$ is $(s/2)$-logarithmically stable.
			\item \label{item:o(1) stable for some s > 0}
			There exists $s_0 > 0$ such that $\{ \semigroup[\gen[\dampco][s_0]] \}_{t \geq 0}$ is $o(1)$ stable.
		\end{eqenumerate}
	\end{quotetheorem}
	Note that $\eqref{item:log stable for any s > 0} \implies \eqref{item:o(1) stable for some s > 0}$ is trivial.
	Also, $\eqref{item:o(1) stable for some s > 0} \implies \eqref{item:d-GCC}$ easily follows from Theorem \ref{theorem:L-S}. 
	\begin{proofof}{$\eqref{item:o(1) stable for some s > 0} \implies \eqref{item:d-GCC}$}
		Let $s_0 > 0$ and suppose that $\{ \semigroup[\gen[\dampco][s_0]] \}_{t \geq 0}$ is $o(1)$ stable.
		Then, by Theorem \ref{theorem:o(1) stable = pointwise strong annihilation}, 
		there exist $\levelofdampco, \levelofm > 0$
		such that $( \setS , \setSigma[1][s_0] )$ is strongly annihilating.
		Here notice that $\measure{ \setSigma[1][s_0] } > 0$. 
		Therefore, using $\eqref{item:L-S SA for some non-null set} \implies \eqref{item:L-S d-GCC}$ of Theorem \ref{theorem:L-S}, 
		we conclude that $\R^d \setminus \setS$ satisfies $d$-GCC.
	\end{proofof}
	We prove the remaining $\eqref{item:d-GCC} \implies \eqref{item:log stable for any s > 0}$ using the following Theorem \ref{theorem:L-S with explicit constant}.
	\begin{quotetheorem}[\cite{MR1840110}*{Theorem 4}] \label{theorem:L-S with explicit constant}
		Let $n \in \N$, $r, \sigma > 0$ and $S, \Sigma \subset \R^d$. 
		If 
		\begin{equation}
		\rho \coloneqq \inf_{a \in \R^d} \frac{ \measure{ B(a, r) \cap ( \R^d \setminus S ) } }{ \measure{ B(a, r) } } > 0 
		\end{equation}
		holds
		and $\Sigma \subset \R^d$ is covered by $n$ balls of radius $\sigma$,
		then 
		$(S, \Sigma)$ is strongly annihilating 
		and
		\begin{equation}
		\anihiconst[S][\Sigma] \leq \mleft( \frac{C^d}{\rho} \mright)^{ d r \sigma \mleft( \frac{C^d}{\rho} \mright)^n - n + \frac{1}{2} } ,
		\end{equation}
	where $C > 0$ denotes some universal constant.
	\end{quotetheorem} 
	\begin{proofof}{$\eqref{item:d-GCC} \implies \eqref{item:log stable for any s > 0}$}
		Let $\levelofdampco > 0$ be such that $\R^d \setminus \setS$ satisfies $d$-GCC.
		Then, since 
		\begin{equation}
			\setSigma[\lambda][2][1] = \set{ \xi \in \R^d }{ \abs{ (\abs{\xi}^2 + 1)^{1/2} - \lambda } < 1 } \subset B( 0, \lambda + 1 ) ,
		\end{equation}
		Theorem \ref{theorem:L-S with explicit constant} implies that there exists $C > 0$ such that $(\setS, \setSigma[\lambda][2][1])$ is strongly annihilating and
		\begin{equation}
			\anihiconst[\setS][\setSigma[\lambda][2][1]] \leq C \exp{( C (\lambda + 1) )}
		\end{equation}
		holds for any $\lambda \geq 0$.
		Therefore, by Theorem \ref{theorem:log stable = exp strong annihilation}, 
		$\{ \semigroup[\gen[\dampco][2]] \}_{t \geq 0}$ is $1$-logarithmically stable.
		Hence, by \eqref{item:log stable for some s then} of Theorem \ref{theorem:stable for some s then for any s}, 
		we conclude that $\{ \semigroup[\gen[\dampco]] \}_{t \geq 0}$ is $(s/2)$-logarithmically stable for any $s > 0$.
	\end{proofof}
	\begin{remark}
		By Theorem \ref{theorem:d-GCC log o(1)}, 
		in the case of $\dampop \colon f \mapsto \sqrt{\dampco} f$, 
		the existence of $s_0 > 0$ such that $\{ \semigroup[\gen[\dampco][s_0]] \}_{t \geq 0}$ is $o(1)$ stable 
		implies 
		the $(s/2)$-logarithmic stability of $\{ \semigroup[\gen[\dampco]] \}_{t \geq 0}$ for any $s > 0$.
		On the other hand, 
		in the case of general bounded operators $\dampop \colon \hilbert \to \hilbert$,
		the existence of $s_0 > 0$ such that $\{ \semigroup[\gen[\dampop][s_0]] \}_{t \geq 0}$ is $o(1)$ stable 
		implies 
		a weaker result, 
		the $o(1)$ stability of $\{ \semigroup \}_{t \geq 0}$ for any $s > 0$, 
		by \eqref{item:o(1) stable for some s then} of Theorem \ref{theorem:stable for some s then for any s}.
		As of this writing, 
		it remains open that 
		the $o(1)$ stability implies the logarithmic stability or not 
		in the case of general bounded operators.
	\end{remark}
	Note that in the case $d = 1$, we have a much stronger result.
	\begin{quotetheorem}[\cite{MR4143391}*{Theorem 1}] \label{theorem:d=1 1-GCC exp poly o(1)}
		Let $d = 1$ and $0 \leq \dampco \in \Linfty(\R)$.
		Then the following are equivalent:
		\begin{eqenumerate}
			\item[\eqref{item:d-GCC}]
			There exists $\levelofdampco > 0$ such that $\R \setminus \setS$ satisfies $1$-GCC.
			\item \label{item:exp stable for any s geq 2 and poly stable for 0 < s < 2}
			$\{ \semigroup[\gen[\dampco]] \}_{t \geq 0}$ is $s/(4 - 2s)$-polynomially stable for any $0 < s < 2$ and exponentially stable for any $s \geq 2$.
			\item[\eqref{item:o(1) stable for some s > 0}]
			There exists $s_0 > 0$ such that $\{ \semigroup[\gen[\dampco][s_0]] \}_{t \geq 0}$ is $o(1)$ stable.
		\end{eqenumerate}
	\end{quotetheorem}
	\begin{proofof}{Theorem \ref{theorem:d=1 1-GCC exp poly o(1)}}
		It is enough to show that $\eqref{item:d-GCC} \implies \eqref{item:exp stable for any s geq 2 and poly stable for 0 < s < 2}$.
		Let $\levelofdampco > 0$ be such that $\R \setminus \setS$ satisfies $1$-GCC.
		Then, since 
		\begin{equation}
			\setSigma[\lambda][2][1] = \set{ \xi \in \R }{ \abs{ (\abs{\xi}^2 + 1)^{1/2} - \lambda } < 1 } \subset B( - \lambda, 2 ) \cup B( \lambda, 2 ) ,
		\end{equation}
		Theorem \ref{theorem:L-S with explicit constant} implies that $\{ (\setS, \setSigma[\lambda][2][1]) \}_{\lambda \geq 0}$ is uniformly strongly annihilating.
		Therefore, by Theorem \ref{theorem:exp stable = uniform strong annihilation} and \eqref{item:exp stable for some s then} of Theorem \ref{theorem:stable for some s then for any s}, 
		we conclude that $\{ \semigroup[\gen[\dampco]] \}_{t \geq 0}$ is 
		exponentially stable for any $s \geq 2$ 
		and $s/(4 - 2s)$-polynomially stable for any $0 < s < 2$.
	\end{proofof}
	The essential difference here is that an $1$-dimensional annulus is just a union of two intervals.
	Recently, \cite{MR4363753} established Theorem \ref{theorem:L-S 1-GCC}, \ref{theorem:L-S with explicit constant annulus ver} and \ref{theorem:d>1 1-GCC exp poly}, 
	analogies of Theorem \ref{theorem:L-S}, \ref{theorem:L-S with explicit constant} and \ref{theorem:d=1 1-GCC exp poly o(1)}, respectively.
	\begin{quotetheorem}[\cite{MR4363753}*{Proposition 6}] \label{theorem:L-S 1-GCC}
		Let $S \subset \R^d$ be open.
		Then the following are equivalent:
		\begin{eqenumerate}
			\item 
			$\R^d \setminus S$ satisfies $1$-GCC.
			\item 
			For any $\sigma > 0$ and $R \in \SO(d)$, $(S, R( [0, \sigma] \times \R^{d - 1} ) )$ is strongly annihilating. 
			\item 
			There exists $\sigma > 0$ such that for any $R \in \SO(d)$, $(S, R( [0, \sigma] \times \R^{d - 1} ) )$ is strongly annihilating. 
		\end{eqenumerate}
		Here $R( [0, \sigma] \times \R^{d - 1} )$ denotes 
		\begin{equation}
			R( [0, \sigma] \times \R^{d - 1} ) \coloneqq \set{ Rx \in \R^d }{ x \in [0, \sigma] \times \R^{d - 1} } .
		\end{equation}
	\end{quotetheorem}
	\begin{quotetheorem}[\cite{MR4363753}*{Theorem 1}] \label{theorem:L-S with explicit constant annulus ver}
		Let $\lambda \geq 0$, $\ell, \sigma > 0$ and $S, \Sigma \subset \R^d$.
		If
		\begin{equation}
		\rho \coloneqq 	\inf_{\substack{ a \in \R^d ,  e \in S^{d-1}}} \frac{ \mathcal{H}^{1}( L(a, e, \ell) \cap (\R^d \setminus S) ) }{ \ell } > 0 
		\end{equation}
		and 
		\begin{equation}
		\Sigma \subset \set{ \xi \in \R^d }{ \abs{ \abs{ \xi } - \lambda } < \sigma } 
		\end{equation} 
		hold, then 
		$(S_\delta, \Sigma)$ is strongly annihilating and
		\begin{equation}
		\anihiconst[S_\delta][\Sigma] \leq C_d \mleft( \frac{\ell}{\delta} \mright)^{d + 1} \mleft( \frac{C_d}{\rho} \mright)^{ C_d \ell \sigma } 
		\end{equation} 
		for any $\delta > 0$, 
		where $C_d > 0$ denotes a constant depending only on $d$ and
		$S_\delta \subset S$ is defined by
		\begin{equation}
		S_\delta \coloneqq \bigcap_{x \in \R^d \setminus S} (\R^d \setminus B(x, \delta) ) . 
		\end{equation} 
	\end{quotetheorem}
	\begin{quotetheorem}[\cite{MR4363753}*{Theorem 4}] \label{theorem:d>1 1-GCC exp poly}
		Let $d \geq 2$ and $0 \leq \dampco \in \Linfty$.
		Then we have the following:
		\begin{eqenumerate}
			\item \label{item:1-GCC is sufficient}
			If there exist $\levelofdampco, \delta > 0$ and $S \subset \R^d$ 
			such that $\R^d \setminus S$ satisfies $1$-GCC and 
			$\setS \subset S_\delta$ holds, 
			then $\{ \semigroup[\gen[\dampco]] \}_{t \geq 0}$ is exponentially stable for any $s \geq 2$ 
			and $s/(4 - 2s)$-polynomially stable for any $0 < s < 2$.
			\item \label{item:1-GCC is necessary}
			If $\{ \semigroup[\gen[\dampco][2]] \}_{t \geq 0}$ is exponentially stable and $\dampco$ is continuous, 
			then there exists $\levelofdampco > 0$ 
			such that $\R^d \setminus \setS$ satisfies $1$-GCC.
		\end{eqenumerate}
	\end{quotetheorem}
	We prove Theorem \ref{theorem:d>1 1-GCC exp poly} 
	using Theorem \ref{theorem:stable for some s then for any s}, \ref{theorem:exp stable = uniform strong annihilation}, \ref{theorem:L-S 1-GCC} and \ref{theorem:L-S with explicit constant annulus ver}.
	\begin{proofof}{Theorem \ref{theorem:d>1 1-GCC exp poly}}
		At first we prove \eqref{item:1-GCC is sufficient}.
		Let $\levelofdampco, \delta > 0$ and $S \subset \R^d$ be
		such that $\R^d \setminus S$ satisfies $1$-GCC and 
		$\setS \subset S_\delta$ holds.
		Then, since 
		\begin{equation}
			\setSigma[\lambda][2][1] 
			= \set{ \xi \in \R^d }{ \abs{ (\abs{\xi}^2 + 1)^{1/2} - \lambda } < 1 } 
			\subset \set{ \xi \in \R^d }{ \abs{ \abs{\xi} - \lambda } < 2 }  ,
		\end{equation}
		Theorem \ref{theorem:L-S with explicit constant annulus ver} implies that $\{ (\setS, \setSigma[\lambda][2][1]) \}_{\lambda \geq 0}$ is uniformly strongly annihilating.
		Therefore, by Theorem \ref{theorem:stable for some s then for any s} and \ref{theorem:exp stable = uniform strong annihilation}, 
		we conclude that $\{ \semigroup[\gen[\dampco]] \}_{t \geq 0}$ is 
		exponentially stable for any $s \geq 2$ 
		and $s/(4 - 2s)$-polynomially stable for any $0 < s < 2$.
		
		Next we prove \eqref{item:1-GCC is necessary}.
		Suppose that $\{ \semigroup[\gen[\dampco][2]] \}_{t \geq 0}$ is exponentially stable. 
		Then, by Theorem \ref{theorem:exp stable = uniform strong annihilation}, 
		there exist $\levelofdampco, \levelofm > 0$ 
		such that $\{ ( \setS , \setSigma[\lambda][2] ) \}_{\lambda \geq 0}$ is uniformly strongly annihilating.
		Since we have
		\begin{align}
			\setSigma[\lambda][2] 
			&= \set{ \xi \in \R^d }{ \abs{ ( \abs{\xi}^2 + 1 )^{1/2} - \lambda } < \levelofm } \\
			&= \set{ \xi \in \R^d }{ ( (\lambda - \levelofm)^{2} - 1 )^{1/2} < \abs{\xi} < ( (\lambda + \levelofm)^{2} - 1 )^{1/2} } 
		\end{align}
		for any $\lambda \geq 1 + \levelofm$
		and 
		\begin{equation}
			\lim_{\lambda \to \infty} ( ( (\lambda + \levelofm)^{2} - 1 )^{1/2} - ( (\lambda - \levelofm)^{2} - 1 )^{1/2} ) = 2 \levelofm ,
		\end{equation}
		a similar argument as that of Proposition \ref{prop:setSigma 0 < s < 2} implies the following property:
		for any bounded sets $\Omega \subset \R^{d - 1}$ and $R \in \SO(d)$, there exist $a \in \R^d$ and $\lambda \geq 0$
		satisfying $a + R( [0, \levelofm] \times \Omega ) \subset \setSigma[\lambda][2]$.
		Now fix $R \in \SO(d)$ and $f \in \hilbert$ 
		such that $\supp \Fourier[f]$ is compact and $\supp \Fourier[f] \subset R( [0, \levelofm] \times \R^{d-1} )$.
		Then we can take $a \in \R^d$ and $\lambda \geq 0$
		such that $\supp \Fourier[\modulation{a} f] = a + \supp \Fourier[f] \subset \setSigma[\lambda][2]$.
		Therefore, since $\{ ( \setS , \setSigma[\lambda][2] ) \}_{\lambda \geq 0}$ is uniformly strongly annihilating, 
		we have
		\begin{equation}
			\norm{f}_{\hilbert} 
			= \norm{ \modulation{a} f }_{\hilbert} 
			\leq C \norm{ \modulation{a} f }_{ L^2( \R^d \setminus \setS ) } 
			= C \norm{ f }_{ L^2( \R^d \setminus \setS ) } , 
		\end{equation}
		where $C = \sup_{\lambda \geq 0} \anihiconst[ \setS ][ \setSigma[\lambda][2] ]$.
		Thus $( \setS, R( [0, \levelofm] \times \R^{d - 1} ) )$ is strongly annihilating for any $R \in \SO(d)$.
		Now note that $\setS$ is open since $\dampco$ is continuous.
		Hence, by Theorem \ref{theorem:L-S 1-GCC}, we conclude that $\R^d \setminus \setS$ satisfies $1$-GCC.
	\end{proofof}
	Finally, we remark that if $0 \leq \dampco \in \Linfty$ is uniformly continuous, 
	then Theorem \ref{theorem:d>1 1-GCC exp poly} is simplified as follows:
	\begin{quotetheorem} \label{theorem:d>1 1-GCC exp poly with uniform continuity}
		Let $d \geq 2$ and $0 \leq \dampco \in \Linfty$ be uniformly continuous.
		Then the following are equivalent:
		\begin{eqenumerate}
			\item[{\eqref{item:1-GCC}}] 
			There exists $\levelofdampco > 0$ such that $\R^d \setminus \setS$ satisfies $1$-GCC.
			\item[{\eqref{item:exp stable for any s geq 2 and poly stable for 0 < s < 2}}]
			$\{ \semigroup[\gen[\dampco]] \}_{t \geq 0}$ is $s/(4 - 2s)$-polynomially stable for any $0 < s < 2$ and exponentially stable for any $s \geq 2$.
			\item
			$\{ \semigroup[\gen[\dampco][2]] \}_{t \geq 0}$ is exponentially stable.
		\end{eqenumerate}
	\end{quotetheorem}
	To see this, let $\levelofdampco > 0$ be such that $\R^d \setminus \setS$ satisfies $1$-GCC.
	Then, by the uniform continuity, there exists $\delta > 0$ such that $\setS[\dampco][\levelofdampco / 2] \subset (\setS)_\delta$.
	Therefore, $\eqref{item:1-GCC} \implies \eqref{item:exp stable for any s geq 2 and poly stable for 0 < s < 2}$ follows from \eqref{item:1-GCC is sufficient} of Theorem \ref{theorem:d>1 1-GCC exp poly}.
	
	\section*{Acknowledgments}
	I am thankful to my colleague Kotaro Inami for introducing me to this problem.

\end{document}